\documentclass[a4paper,12pt]{article}

\usepackage{amsmath}
\usepackage{ntheorem}
\usepackage{amssymb}
\usepackage{latexsym}
\usepackage[all]{xy}\UseComputerModernTips
\usepackage{graphicx}

\theoremstyle{plain}
\newtheorem{prp}{Proposition}[section]
\newtheorem{thm}[prp]{Theorem}
\newtheorem{lem}[prp]{Lemma}
\newtheorem{cor}[prp]{Corollary}

\theoremstyle{plain}
\theorembodyfont{\normalfont\rm}
\newtheorem{dfn}[prp]{Definition}
\newtheorem{rem}[prp]{Remark}
\newtheorem{exa}[prp]{Example}

\theoremstyle{nonumberplain}
\theoremheaderfont{\normalfont\itshape}
\theoremseparator{.}
\theorembodyfont{\normalfont}
\newtheorem{proof}{Proof}

\newcommand{\qed}{\hfill $\Box$}

\makeatletter
\renewcommand{\theequation}{%
\thesection.\arabic{equation}}
\@addtoreset{equation}{section}
\makeatother

\newcommand{\ZZ}{{\mathbb Z}}
\newcommand{\NN}{{\mathbb N}}
\newcommand{\RR}{{\mathbb R}}
\newcommand{\CC}{{\mathbb C}}
\newcommand{\PP}{{\mathbb P}}

\newcommand{\LL}{{\mathbb L}}

\renewcommand{\d}{{\rm dim}}

\newcommand{\Vol}{{\rm Vol}}

\newcommand{\e}{\varepsilon}

\newcommand{\id}{{\rm id}}
\newcommand{\supp}{{\rm supp}}

\newcommand{\Int}{{\rm Int}}

\newcommand{\relint}{{\rm rel.int}}

\newcommand{\Cone}{{\rm Cone}}
\newcommand{\bif}{{\rm bif}}
\renewcommand{\phi}{{\varphi}}

\newcommand{\Db}{{\bf D}^{b}}
\newcommand{\Dbc}{{\bf D}_{c}^{b}}
\newcommand{\Kbc}{{\bf K}_{c}^{b}}

\newcommand{\F}{{\cal F}}
\newcommand{\G}{{\cal G}}

\newcommand{\CF}{{\rm CF}}
\newcommand{\J}{{S}}

\newcommand{\RG}{R\varGamma}

\newcommand{\tl}[1]{\widetilde{#1}}

\newcommand{\simto}{\overset{\sim}{\longrightarrow}}

\newcommand{\dsum}{\displaystyle \sum}
\newcommand{\dint}{\displaystyle \int}

\renewcommand{\(}{\left(}
\renewcommand{\)}{\right)}

\newcommand{\longhookrightarrow}{\DOTSB\lhook\joinrel\longrightarrow}
\newcommand{\longtwoheadrightarrow}{\relbar\joinrel\twoheadrightarrow}

\title{Monodromy zeta functions at infinity, Newton polyhedra and constructible sheaves }
\author{Yutaka \textsc{Matsui}\footnote{Department of Mathematics, Kinki University, 3-4-1, Kowakae, Higashi-Osaka, Osaka, 577-8502, Japan.} \and Kiyoshi \textsc{Takeuchi}\footnote{Institute of Mathematics, University  of Tsukuba, 1-1-1, Tennodai, Tsukuba, Ibaraki, 305-8571, Japan.}}

\date{}

\begin{document}

\maketitle

\begin{abstract}
By using sheaf-theoretical methods such as constructible sheaves, we generalize the formula of Libgober-Sperber \cite{L-S} concerning the zeta functions of monodromy at infinity of polynomial maps into various directions. In particular, some formulas for the zeta functions of global monodromy along the fibers of bifurcation points of polynomial maps will be obtained.
\end{abstract}

\section{Introduction}\label{sec:1}

After two fundamental papers \cite{Broughton} and \cite{S-T-1}, numerous papers have been written to study the global behavior of polynomial maps $f \colon \CC^n \longrightarrow \CC$. For a polynomial map $f \colon \CC^n \longrightarrow \CC$, it is well-known that there exists a finite subset $B \subset \CC$ of $\CC$ such that the restriction
\begin{equation}
\CC^n \setminus f^{-1}(B) \longrightarrow \CC \setminus B
\end{equation}
of $f$ is a locally trivial fibration over $\CC \setminus B$. We denote by $B_f$ the smallest subset $B \subset \CC$ satisfying this condition and call its elements bifurcation points of $f$. Let $F=f^{-1}(c)$ ($c\in \CC \setminus B_f$) be the generic fiber. Then we obtain a monodromy representations
\begin{equation}
\Psi_j \colon \pi_1(\CC \setminus B_f; c) \longrightarrow {\rm Aut}(H_j(F; \CC))
\end{equation}
for $j \geq 0$. Despite important contributions by many mathematicians, these monodromy representations still remain mysterious. The aim of this paper is to prove some formulas for the zeta functions associated to these representations by sheaf-theoretical methods. Our results will be completely described by certain Newton polyhedra associated to $f$. In order to explain our results more precisely, first let us recall the definition of the monodromy at infinity of $f$ studied by many mathematicians. Let $C_R=\{x\in \CC\ |\ |x|=R\}$ ($R\gg 0$) be a sufficiently large circle in $\CC$ such that $B_f\subset \{x \in \CC\ |\ |x|<R\}$. Then by restricting the locally trivial fibration 
\begin{equation}
\CC^n \setminus f^{-1}(B_f) \longrightarrow \CC \setminus B_f
\end{equation}
to $C_R$ we obtain a geometric monodromy automorphism 
\begin{equation}
\Phi_f^{\infty} \colon f^{-1}(R) \simto f^{-1}(R).
\end{equation}
We define the monodromy zeta function at infinity $\zeta_{f}^{\infty}(t) \in \CC(t)^*$ of $f$ by 
\begin{equation}
\zeta_{f}^{\infty}(t):=\prod_{j=0}^{\infty} \det(\id -t\Phi_j)^{(-1)^j},
\end{equation}
where 
\begin{equation}
\Phi_j \colon H^j(f^{-1}(R) ;\CC) \overset{\sim}{\longrightarrow} H^j(f^{-1}(R) ;\CC) \ \ (j=0,1,\ldots)
\end{equation}
are the isomorphisms induced by $\Phi_f^{\infty}$. Various formulas for monodromy zeta functions at infinity were obtained by Libgober-Sperber \cite{L-S}, Garc{\'i}a-L{\'o}pez-N{\'e}methi \cite{L-N}, Gusein-Zade-Luengo-Melle-Hern{\'a}ndez \cite{GLM2}, \cite{GLM1} and Siersma-Tib{\u a}r \cite{S-T-2} etc. In particular, Libgober-Sperber \cite{L-S} proved a beautiful formula which expresses the monodromy zeta function at infinity of a polynomial $f(x)\in \CC [x_1, x_2, \ldots , x_n]$ in terms of its Newton polytope at infinity. In this paper, by using some functorial properties of nearby cycle and constructible sheaves, we first give a new proof to this Libgober-Sperber's formula and generalize it to non-convenient polynomials (see Theorem \ref{thm:3-5} (i)). We eliminate the points of indeterminacy of the meromorphic extension $\tl{f}$ of $f(x)\in \CC [x_1, x_2, \ldots , x_n]$ to the compactification of $(\CC^*)^n$ used in the proof of \cite{L-S} by blowing up it several times and prove a formula for the monodromy zeta function $\zeta_f^{\infty}(t)$ by calculating that of a constructible sheaf on the resulting complex manifold. This functorial proof leads us to various applications. For example in Section \ref{sec:3} and \ref{sec:4}, for any bifurcation point $b \in B_f$ of $f(x) \in \CC[x_1,x_2,\ldots,x_n]$ we prove explicit formulas for the global monodromy zeta function $\zeta_f^b(t) \in \CC (t)^*$ along the fiber $f^{-1}(b)$ obtained by restricting the locally trivial fibration $\CC^n \setminus f^{-1}(B_f) \longtwoheadrightarrow \CC\setminus B_f$ to a small loop $\{x\in \CC\ |\ |x-b |= \varepsilon \}$ ($0<\varepsilon \ll 1$) around $b \in B_f\subset \CC$. One of the most unexpected results we obtain in Section \ref{sec:3} is that the constant term $a=a_0 \in \CC$ of a non-convenient polynomial $f(x)= \sum_{v\in \ZZ_{\geq 0}^n} a_vx^v$ is a bifurcation point of $f$ in general. In Corollary \ref{cor:4-6} and \ref{cor:4-10}, we will give some formulas which express the jumping number
\begin{equation}
\chi (f^{-1}(a))-\chi (f^{-1}(c)) \in \ZZ
\end{equation}
of the Euler characteristic of the central fiber $f^{-1}(a)$ from that of the general fiber $f^{-1}(c)$ in terms of a Newton polyhedron of $f-a$. Note that Siersma-Tib{\u a}r \cite{S-T-2} also proved formulas for the global monodromy zeta function $\zeta_f^b(t)$ along the fiber $f^{-1}(b)$ ($b \in B_f$) by completely different methods. However the results in \cite{S-T-2} are described by the polar curves associated with $f$ and their relation with the Newton polyhedron of $f$ is not clear. Note also that in \cite[Theorem 4]{GLM1} Gusein-Zade-Luengo-Melle-Hern{\'a}ndez obtained a formula for $\zeta_f^b(t)$ by using the meromorphic extension $\tl{f-b}$ of $f-b$ to the projective space $\PP^n$. In order to use their formula, we have to compute the monodromy zeta functions of the meromorphic function $\tl{f-b}$ at all points of $\overline{\{f=b\}} \cap H_{\infty}$, where $H_{\infty}=\PP^n \setminus \CC^n$ is the hyperplane at infinity. Our formulas for $\zeta_f^b(t)$ in Theorem \ref{thm:4-5}, \ref{thm:4-9} and \ref{thm:4-12} are directly described by the Newton polyhedron of $f$. In the final section, we prove also a generalization of Libgober-Sperber's formula to polynomial maps $f=(f_1,f_2, \ldots , f_k)\colon \CC^n \longrightarrow \CC^k$ ($1 \leq k \leq n$). Namely we study the global Milnor fiber $f^{-1}(c)$ of $f$ associated with the complete intersection subvariety $\{f_1=f_2= \cdots =f_k=0\}$ in $\CC^n$, where $c \in \CC^k$ is a generic point, i.e. a point outside the discriminant $D \subset \CC^k$ of $f \colon \CC^n \longrightarrow \CC^k$. The results we obtain in Section \ref{sec:5} are described also by certain Newton polyhedra defined by $f_1, f_2, \ldots , f_k \in \CC[x_1,x_2,\ldots,x_n]$ and can be considered as the global versions of the theorem of Kirillov \cite{Kirillov} and Oka \cite{Oka-2}. Finally, let us mention that the sheaf-theoretical methods we used in this paper can be applied also to other problems. For example, in \cite{M-T-new2} we computed the monodromy zeta functions of Milnor fibers over general singular toric varieties. In another recent paper \cite{M-T-new1}, some applications of our methods to $A$-discriminant varieties are given. Moreover in \cite{Takeuchi-2} we used the methods developed in this paper to obtain a formula for the monodromy at infinity of $A$-hypergeometric functions.

\section{Preliminary notions and results}\label{sec:2}

In this section, we introduce basic notions and results which will be used in this paper. In this paper, we essentially follow the terminology of \cite{Dimca} and \cite{K-S}. For example, for a topological space $X$ we denote by $\Db(X)$ the derived category whose objects are bounded complexes of sheaves of $\CC_X$-modules on $X$.

\begin{dfn}\label{dfn:2-1}
Let $X$ be an algebraic variety over $\CC$. Then
\begin{enumerate}
\item We say that a sheaf $\F$ on $X$ is constructible if there exists a stratification $X=\bigsqcup_{\alpha} X_{\alpha}$ of $X$ such that $\F|_{X_{\alpha}}$ is a locally constant sheaf of finite rank for any $\alpha$.
\item We say that an object $\F$ of $\Db(X)$ is constructible if the cohomology sheaf $H^j(\F)$ of $\F$ is constructible for any $j \in \ZZ$. We denote by $\Dbc(X)$ the full subcategory of $\Db(X)$ consisting of constructible objects $\F$.
\end{enumerate}
\end{dfn}

Recall that for any morphism $f \colon X \longrightarrow Y$ of algebraic varieties over $\CC$ there exists a functor
\begin{equation}
Rf_* \colon \Db(X) \longrightarrow \Db(Y)
\end{equation}
of direct images. This functor preserves the constructibility and we obtain also a functor
\begin{equation}
Rf_* \colon \Dbc(X) \longrightarrow \Dbc(Y).
\end{equation}
For other basic operations $Rf_!$, $f^{-1}$, $f^!$ etc. in derived categories, see \cite{K-S} for the detail.

Next we introduce the notion of constructible functions and explain its relation with that of constructible sheaves.

\begin{dfn}\label{dfn:2-2}
Let $X$ be an algebraic variety over $\CC$ and $G$ an abelian group. Then we say a $G$-valued function $\rho \colon X \longrightarrow G$ on $X$ is constructible if there exists a stratification $X=\bigsqcup_{\alpha} X_{\alpha}$ of $X$ such that $\rho|_{X_{\alpha}}$ is constant for any $\alpha$. We denote by $\CF_G(X)$ the abelian group of $G$-valued constructible functions on $X$.
\end{dfn}

Let $\CC(t)^*=\CC(t) \setminus \{0\}$ be the multiplicative group of the function field $\CC(t)$ of the scheme $\CC$. In this paper, we consider $\CF_G(X)$ only for $G=\ZZ$ or $\CC(t)^*$. For a $G$-valued constructible function $\rho \colon X \longrightarrow G$, we take a stratification $X=\bigsqcup_{\alpha}X_{\alpha}$ of $X$ such that $\rho|_{X_{\alpha}}$ is constant for any $\alpha$ as above. Denoting the Euler characteristic of $X_{\alpha}$ by $\chi(X_{\alpha})$, we set
\begin{equation}
\int_X \rho :=\dsum_{\alpha}\chi(X_{\alpha}) \cdot \rho(x_{\alpha}) \in G,
\end{equation}
where $x_{\alpha}$ is a reference point in $X_{\alpha}$. Then we can easily show that $\int_X\rho \in G$ does not depend on the choice of the stratification $X=\bigsqcup_{\alpha} X_{\alpha}$ of $X$. Hence we obtain a homomorphism
\begin{equation}
\int_X \colon \CF_G(X) \longrightarrow G
\end{equation}
of abelian groups. For $\rho \in \CF_G(X)$, we call $\int_X \rho \in G$ the topological (Euler) integral of $\rho$ over $X$. More generally, for any morphism $f \colon X \longrightarrow Y$ of algebraic varieties over $\CC$ and $\rho \in \CF_G(X)$, we define the push-forward $\int_f \rho \in \CF_G(Y)$ of $\rho$ by
\begin{equation}
\( \int_f \rho \) (y):=\int_{f^{-1}(y)} \rho
\end{equation}
for $y \in Y$. This defines a homomorphism
\begin{equation}
\int_f \colon \CF_G(X) \longrightarrow \CF_G(Y)
\end{equation}
of abelian groups. If $G=\ZZ$, these operations $\int_X$ and $\int_f$ correspond to the ones $\RG(X;\ \cdot \ )$ and $Rf_*$ respectively in the derived categories as follows. For an algebraic variety $X$ over $\CC$, consider a free abelian group
\begin{equation}
\ZZ(\Dbc(X)):=\left\{ \left. \dsum_{j \colon \text{finite}}a_j [\F_j] \ \right| \ a_j \in \ZZ, \ \F_j \in \Dbc(X)\right\}
\end{equation}
generated by the objects $\F_j \in \Dbc(X)$ in $\Dbc(X)$ and take its subgroup
\begin{equation}
R:=\langle [\F_2]-[\F_1]-[\F_3] | \text{$\F_1 \longrightarrow \F_2 \longrightarrow \F_3 \overset{+1}{\longrightarrow}$ is a distinguished triangle} \rangle \subset \ZZ(\Dbc(X)).
\end{equation}
We set $\Kbc(X):=\ZZ(\Dbc(X))/R$ and call it the Grothendieck group of $\Dbc(X)$. Then the following result is well-known (see for example \cite[Theorem 9.7.1]{K-S}).

\begin{thm}\label{thm:2-3}
The homomorphism
\begin{equation}
\chi_X \colon \Kbc(X) \longrightarrow \CF_{\ZZ}(X)
\end{equation}
defined by taking the local Euler-Poincar{\'e} indices:
\begin{equation}
\chi_X([\F])(x):=\dsum_{j \in \ZZ} (-1)^j \d_{\CC}H^j(\F)_x \hspace{5mm}(x 
\in X)
\end{equation}
is an isomorphism.
\end{thm}

For any morphism $f \colon X \longrightarrow Y$ of algebraic varieties over $\CC$, there exists also a commutative diagram
\begin{equation}
\xymatrix{
\Kbc(X) \ar[r]^{Rf_*} \ar[d]^{\wr}_{\chi_X}& \Kbc(Y) \ar[d]^{\wr}_{\chi_Y} \\
\CF_{\ZZ}(X) \ar[r]^{\int_f} & \CF_{\ZZ}(Y).}
\end{equation}
In particular, if $Y$ is the one-point variety $\{{\rm pt}\}$ ($\Kbc(Y) \simeq \CF_{\ZZ}(Y) \simeq \ZZ$), we obtain a commutative diagram
\begin{equation}
\xymatrix@R=2.5mm@C=30mm{
\Kbc(X) \ar[rd]^{\chi(\RG(X; \ \cdot \ ))} \ar[dd]_{\wr}^{\chi_X}& \\
 & \ZZ. \\
\CF_{\ZZ}(X) \ar[ru]^{\int_X}}
\end{equation}

Among various operations in derived categories, the following nearby cycle functor introduced by Deligne will be frequently used in this paper (see \cite[Section 4.2]{Dimca} for an excellent survey of this subject).

\begin{dfn}\label{dfn:2-4}
Let $f \colon X \longrightarrow \CC$ be a non-constant regular function on an algebraic variety $X$ over $\CC$. Set $X_0:= \{x\in X\ |\ f(x)=0\} \subset X$ and let $i_X \colon X_0 \longhookrightarrow X$, $j_X \colon X \setminus X_0 \longhookrightarrow X$ be inclusions. Let $p \colon \tl{\CC^*} \longrightarrow \CC^*$ be the universal covering of $\CC^* =\CC \setminus \{0\}$ ($\tl{\CC^*} \simeq \CC$) and consider the Cartesian square
\begin{equation}\label{eq:2-13}
\xymatrix@R=2.5mm@C=2.5mm{
\tl{X \setminus X_0} \ar[rr] \ar[dd]^{p_X} & &\tl{\CC^*} \ar[dd]^p \\
 & \Box & \\
X \setminus X_0 \ar[rr]^f & & \CC^*.}
\end{equation}
Then for $\F \in \Dbc(X)$ we set
\begin{equation}
\psi_f(\F) := i_X^{-1}R(j_X \circ p_X)_*(j_X \circ p_X)^{-1}\F \in \Db(X_0)
\end{equation}
and call it the nearby cycle of $\F$. 
\end{dfn}

Since the nearby cycle functor preserves the constructibility, in the above situation we obtain a functor
\begin{equation}
\psi_f \colon \Dbc(X) \longrightarrow \Dbc(X_0).
\end{equation}

As we see in the next proposition, the nearby cycle functor $\psi_f$ generalizes the classical notion of Milnor fibers. First, let us recall the definition of Milnor fibers and Milnor monodromies over singular varieties (see for example \cite{Takeuchi} for a review on this subject). Let $X$ be a subvariety of $\CC^m$ and $f \colon X \longrightarrow \CC$ a non-constant regular function on $X$. Namely we assume that there exists a polynomial function $\tl{f} \colon \CC^m \longrightarrow \CC$ on $\CC^m$ such that $\tl{f}|_X=f$. For simplicity, assume also that the origin $0 \in \CC^m$ is contained in $X_0=\{x \in X \ |\ f(x)=0\}$. Then the following lemma is well-known (see for example \cite[Definition 1.4]{Massey}).

\begin{lem}\label{lem:2-5}
For sufficiently small $\e >0$, there exists $\eta_0 >0$ with $0<\eta_0 \ll \e$ such that for $0 < \forall \eta <\eta_0$ the restriction of $f$:
\begin{equation}
X \cap B(0;\e) \cap \tl{f}^{-1}(D(0;\eta) \setminus \{0\}) \longrightarrow D(0;\eta) \setminus \{0\}
\end{equation}
is a topological fiber bundle over the punctured disk $D(0;\eta) \setminus \{0\}:=\{ z \in \CC \ |\ 0<|z|<\eta\}$, where $B(0;\e)$ is the open ball in $\CC^m$ with radius $\e$ centered at the origin.
\end{lem}

\begin{dfn}\label{dfn:2-6}
A fiber of the above fibration is called the Milnor fiber of the function $f \colon X\longrightarrow \CC$ at $0 \in X$ and we denote it by $F_0$.
\end{dfn}

\begin{prp}{\rm \bf(\cite[Proposition 4.2.2]{Dimca})}\label{prp:2-7}
There exists a natural isomorphism
\begin{equation}
H^j(F_0;\CC) \simeq H^j(\psi_f(\CC_X))_0
\end{equation}
for any $j \in \ZZ$.
\end{prp}

By this proposition, we can study the cohomology groups $H^j(F_0;\CC)$ of the Milnor fiber $F_0$ by using sheaf theory. Recall also that in the above situation, as in the same way as the case of polynomial functions over $\CC^n$ (see \cite{Milnor}), we can define the Milnor monodromy operators
\begin{equation}
\Phi_j \colon H^j(F_0;\CC) \overset{\sim}{\longrightarrow} H^j(F_0;\CC) \hspace{5mm}(j=0,1,\ldots)
\end{equation}
and the zeta-function
\begin{equation}
\zeta_{f,0}(t):=\prod_{j=0}^{\infty} \det(\id -t\Phi_j)^{(-1)^j}
\end{equation}
associated with it. Since the above product is in fact finite, $\zeta_{f,0}(t)$ is a rational function of $t$ and its degree in $t$ is the topological Euler characteristic $\chi(F_0)$ of the Milnor fiber $F_0$. Similarly, also for any $y \in X_0=\{x \in X \ |\ f(x)=0\}$ we can define $F_y$ and $\zeta_{f,y}(t) \in \CC (t)^*$. This classical notion of Milnor monodromy zeta functions can be also generalized as follows.

\begin{dfn}\label{dfn:2-8}
Let $f \colon X \longrightarrow \CC$ be a non-constant regular function on $X$ and $\F \in \Dbc(X)$. Set $X_0 :=\{x\in X\ |\ f(x)=0\}$. Then there exists a monodromy automorphism
\begin{equation}
\Phi(\F) \colon \psi_f(\F) \simto \psi_f(\F)
\end{equation}
of $\psi_f(\F)$ in $\Dbc(X_0)$ associated with a generator of the group ${\rm Deck}(\tl{\CC^*}, \CC^*)\simeq \ZZ$ of the deck transformations of $p \colon \tl{\CC^*} \longrightarrow \CC^*$ in the diagram \eqref{eq:2-13}. We define a $\CC(t)^*$-valued constructible function $\zeta_f(\F) \colon X_0 \longrightarrow \CC(t)^* $ on $X_0$ by
\begin{equation}
\zeta_{f,x}(\F)(t):=\prod_{j \in \ZZ} \det\(\id -t\Phi(\F)_{j,x}\)^{(-1)^j}
\end{equation}
for $x \in X_0$, where $\Phi(\F)_{j,x} \colon (H^j(\psi_f(\F)))_x \simto (H^j(\psi_f(\F)))_x$ is the stalk at $x \in X_0$ of the sheaf homomorphism
\begin{equation}
\Phi(\F)_j \colon H^j(\psi_f(\F)) \simto H^j(\psi_f(\F))
\end{equation}
associated with $\Phi(\F)$.
\end{dfn}

The following proposition will play a crucial role in the proof of our main theorems. For the proof, see for example, \cite[p.170-173]{Dimca} and \cite{Schurmann}.

\begin{prp}\label{prp:2-9}
Let $\pi \colon Y \longrightarrow X$ be a proper morphism of algebraic varieties over $\CC$ and $f \colon X \longrightarrow \CC$ a non-constant regular function on $X$. Set $g:=f \circ \pi \colon Y \longrightarrow \CC$, $X_0:=\{x\in X\ |\ f(x)=0\}$ and $Y_0:=\{y\in Y\ |\ g(y)=0\}=\pi^{-1}(X_0)$. Then for any $\G\in \Dbc(Y)$ we have
\begin{equation}
\int_{\pi|_{Y_0}} \zeta_g(\G) =\zeta_f(R\pi_*\G)
\end{equation}
in $\CF_{\CC(t)^*}(X_0)$, where
\begin{equation}
\int_{\pi|_{Y_0}}\colon \CF_{\CC(t)^*}(Y_0) \longrightarrow \CF_{\CC(t)^*}(X_0)
\end{equation}
is the push-forward of $\CC(t)^*$-valued constructible functions by $\pi|_{Y_0} \colon Y_0 \longrightarrow X_0$.
\end{prp}

Finally we recall Bernstein-Khovanskii-Kushnirenko's theorem \cite{Khovanskii}.

\begin{dfn}
Let $g(x)=\sum_{v \in \ZZ^n} a_vx^v$ be a Laurent polynomial on $(\CC^*)^n$ ($a_v\in \CC$). 
\begin{enumerate}
\item We call the convex hull of $\supp(g):=\{v\in \ZZ^n \ |\ a_v\neq 0\} \subset \ZZ^n \subset \RR^n$ in $\RR^n$ the Newton polygon of $g$ and denote it by $NP(g)$.
\item For a vector $u\in \RR^n$, we set
\begin{equation}
\Gamma(g;u):=\left\{ v\in NP(g)\ \left| \ \langle u,v\rangle =\min_{w\in NP(g)} \langle u,w\rangle \right.\right\},
\end{equation}
where for $u=(u_1,\ldots,u_n)$ and $v=(v_1,\ldots, v_n)$ we set $\langle u,v\rangle =\sum_{i=1}^n u_iv_i$.
\item For a vector $u \in \RR^n$, we define the $u$-part of $g$ by
\begin{equation}
g^u(x):=\dsum_{v \in \Gamma(g;u)} a_vx^v.
\end{equation}
\end{enumerate}
\end{dfn}

\begin{dfn}
Let $g_1, g_2, \ldots , g_p$ be Laurent polynomials on $(\CC^*)^n$. Then we say that the subvariety $Z^*=\{ x\in (\CC^*)^n \ |\ g_1(x)=g_2(x)= \cdots =g_p(x)=0 \}$ of $(\CC^*)^n$ is non-degenerate complete intersection if for any covector $u \in \ZZ^n$ the $p$-form $dg_1^u \wedge dg_2^u \wedge \cdots \wedge dg_p^u$ does not vanish on $\{ x\in (\CC^*)^n \ |\ g_1^u(x)= \cdots =g_p^u(x)=0 \}$.
\end{dfn}

\begin{thm}[\cite{Khovanskii}]\label{thm:2-10}
Let $g_1, g_2, \ldots , g_p$ be Laurent polynomials on $(\CC^*)^n$. Assume that the subvariety $Z^*=\{ x\in (\CC^*)^n \ |\ g_1(x)=g_2(x)= \cdots =g_p(x)=0 \}$ of $(\CC^*)^n$ is non-degenerate complete intersection. Set $\Delta_i:=NP(g_i)$ for $i=1,\ldots, p$. Then we have
\begin{equation}
\chi(Z^*)=(-1)^{n-p}\dsum_{\begin{subarray}{c} a_1,\ldots,a_p \geq 1\\ a_1+\cdots +a_p=n\end{subarray}}\Vol_{\ZZ}(\underbrace{\Delta_1,\ldots,\Delta_1}_{\text{$a_1$-times}},\ldots,\underbrace{\Delta_p,\ldots,\Delta_p}_{\text{$a_p$-times}}), 
\end{equation}
where $\Vol_{\ZZ}(\underbrace{\Delta_1,\ldots,\Delta_1}_{\text{$a_1$-times}},\ldots,\underbrace{\Delta_p,\ldots,\Delta_p}_{\text{$a_p$-times}})\in \ZZ$ is the normalized $n$-dimensional mixed volume of $\underbrace{\Delta_1,\ldots,\Delta_1}_{\text{$a_1$-times}},\ldots,\underbrace{\Delta_p,\ldots,\Delta_p}_{\text{$a_p$-times}}$ with respect to the lattice $\ZZ^n \subset \RR^n$.
\end{thm}

\begin{rem}\label{rem:2-13}
Let $Q_1,Q_2,\ldots,Q_n$ be integral polytopes in $(\RR^n, \ZZ^n)$. Then their normalized $n$-dimensional mixed volume $\Vol_{\ZZ}(Q_1,Q_2,\ldots,Q_n) \in \ZZ$ is given by the formula 
\begin{eqnarray}
\lefteqn{n! \Vol_{\ZZ}(Q_1, Q_2, \ldots , Q_n)}\nonumber\\
&=&\Vol_{\ZZ}(Q_1+ Q_2+ \cdots + Q_n)-\sum_{i=1}^n \Vol_{\ZZ}(Q_1+ \cdots +Q_{i-1}+Q_{i+1} + \cdots + Q_n)\nonumber\\
& &+\sum_{1 \leq i<j \leq n} \Vol_{\ZZ}(Q_1+ \cdots +Q_{i-1}+Q_{i+1} + \cdots +Q_{j-1}+Q_{j+1} + \cdots + Q_n)\nonumber\\
& &+\cdots + (-1)^{n-1} \sum_{i=1}^n \Vol_{\ZZ}(Q_i),
\end{eqnarray}
where $\Vol_{\ZZ}(\ \cdot\ )\in \ZZ$ is the normalized $n$-dimensional volume.
\end{rem}

\section{Monodromy at infinity and Newton polyhedra}\label{sec:3}

In this section, we apply our methods to the monodromy at infinity of polynomials on $\CC^n$ studied by Gusein-Zade-Luengo-Melle-Hern{\'a}ndez \cite{GLM2}, \cite{GLM1}, Libgober-Sperber \cite{L-S}, Garc{\'i}a-L{\'o}pez-N{\'e}methi \cite{L-N} and Siersma-Tib{\u a}r \cite{S-T-1}, \cite{S-T-2} etc. Hereafter we denote $\ZZ_{\geq 0}$ by $\ZZ_+$.

\begin{dfn}[\cite{L-S}]\label{dfn:3-1}
Let $f(x)=\sum_{v \in \ZZ_+^n} a_vx^v \in \CC[x_1,x_2,\ldots, x_n]$ ($a_v \in \CC$) be a polynomial on $\CC^n$. We call the convex hull of $\{0\}\cup NP(f)$ in $\RR^n$ the Newton polygon of $f$ at infinity and denote it by $\Gamma_{\infty}(f)$.
\end{dfn}

For a subset $\J\subset \{1,2,\ldots,n\}$ of $\{1,2,\ldots, n\}$, let us set
\begin{equation}
\RR^{\J}:=\{ v=(v_1,v_2,\ldots,v_n)\in \RR^n\ |\ \text{$v_i=0$ for $\forall i \notin \J$}\}.
\end{equation}
We set also $\Gamma_{\infty}^{\J}(f)= \Gamma_{\infty}(f) \cap \RR^{\J}$.

\begin{dfn}\label{dfn:3-2}
We say that a polynomial $f(x)\in \CC[x_1,x_2,\ldots,x_n]$ on $\CC^n$ satisfies the condition ($\ast$) if $\Gamma_{\infty}^{\J}(f)=\{0\}$ or the dimension of $\Gamma_{\infty}^{\J}(f)$ is maximal i.e. equal to $\sharp {\J}$ for any subset $\J$ of $\{1,2,\ldots,n\}$.
\end{dfn}

Recall that a polynomial $f(x)$ on $\CC^n$ is called convenient if the dimension of $\Gamma_{\infty}^{\J}(f)$ is equal to $\sharp {\J}$ for any ${\J} \subset \{1,2,\ldots,n\}$. Therefore convenient polynomials on $\CC^n$ satisfy our condition ($\ast$).

\begin{dfn}[\cite{Kushnirenko}]\label{dfn:3-3}
We say that a polynomial $f(x)=\sum_{v\in \ZZ_+^n} a_vx^v\in \CC[x_1,x_2,\ldots,x_n]$ ($a_v\in \CC$) on $\CC^n$ is non-degenerate at infinity if for any face $\gamma$ of $\Gamma_{\infty}(f)$ such that $0 \notin \gamma$ the complex hypersurface
\begin{equation}
\{x=(x_1,x_2,\ldots,x_n) \in (\CC^*)^n\ |\ f_{\gamma}(x)=0\}
\end{equation}
in $(\CC^*)^n$ is smooth and reduced, where we set $f_{\gamma}(x)=\sum_{v \in \gamma \cap \ZZ_+^n} a_vx^v \in \CC[x_1,x_2,\ldots,x_n]$.
\end{dfn}

Now let $f(x)$ be a polynomial on $\CC^n$. Then it is well-known that there exists a finite subset $B \subset \CC$ of $\CC$ such that the restriction
\begin{equation}
\CC^n \setminus f^{-1}(B) \longrightarrow \CC \setminus B
\end{equation}
of $f$ is a locally trivial fibration. We denote by $B_f$ the smallest subset $B \subset \CC$ verifying this condition and call it the bifurcation set of $f$. Our objective here is the study of the following monodromy zeta functions.

\begin{dfn}\label{dfn:3-4}
\begin{enumerate}
\item Take a sufficiently large circle $C_R=\{x\in \CC\ |\ |x|=R\}$ ($R\gg 0$) in $\CC$ such that $B_f\subset \{x \in \CC\ |\ |x|<R\}$. By restricting the locally trivial fibration $\CC^n \setminus f^{-1}(B_f) \longtwoheadrightarrow \CC\setminus B_f$ to $C_R \subset \CC \setminus B_f$, we obtain the geometric monodromy at infinity
\begin{equation}
\Phi_f^{\infty} \colon f^{-1}(R) \simto f^{-1}(R).
\end{equation}
We denote the zeta function associated with $\Phi_f^{\infty}$ by $\zeta_f^{\infty}(t) \in \CC(t)^*$ and call it the monodromy zeta function at infinity of $f$. 
\item For a bifurcation point $b\in B_f$ of $f$, take a small circle $C_{\e}(b)=\{x\in \CC\ |\ |x-b|=\e\}$ ($0<\e\ll 1$) around $b$ such that $B_f\cap \{x \in \CC\ |\ |x-b|\leq\e\}=\{b\}$. We denote by $\zeta_f^b(t)\in \CC(t)^*$ the zeta function associated with the geometric monodromy
\begin{equation}
\Phi_f^b \colon f^{-1}(b+\e) \simto f^{-1}(b+\e)
\end{equation}
obtained by the restriction of $\CC^n \setminus f^{-1}(B_f) \longtwoheadrightarrow \CC \setminus B_f$ to $C_{\e}(b)\subset \CC\setminus B_f$. We call $\zeta_f^b(t)$ the monodromy zeta function of $f$ along the fiber $f^{-1}(b)$.
\end{enumerate}
\end{dfn}

To compute the monodromy zeta function $\zeta_f^b(t)\in \CC(t)^*$ of $f$ along the fiber $f^{-1}(b)$ of $b\in B_f$, it is very useful to consider first the following rational function $\tl{\zeta_f^b}(t)\in \CC(t)^*$. Let $f^{-1}(b)=\bigsqcup_{\alpha}Z_{\alpha}$ be a stratification of $f^{-1}(b)=\{ f-b=0\}$ such that the local monodromy zeta function $\zeta_{f-b}(t)$ of $f$ is constant on each stratum $Z_{\alpha}$. Denote the value of $\zeta_{f-b}(t)$ on $Z_{\alpha}$ by $\zeta_{\alpha}(t) \in \CC(t)^*$. Then the following definition does not depend on the stratification $f^{-1}(b)=\bigsqcup_{\alpha}Z_{\alpha}$.

\begin{dfn}
We set
\begin{equation}
\tl{\zeta_f^b}(t):= \int_{f^{-1}(b)} \zeta_{f-b}(t)= \prod_{\alpha} \{ \zeta_{\alpha}(t) \}^{\chi(Z_{\alpha})}\in \CC(t)^*
\end{equation}
and call it the finite part of $\zeta_f^b(t)$.
\end{dfn}

As is clear from the definition, the finite part $\tl{\zeta_f^b}(t)$ of $\zeta_f^b(t)$ is calculated only by the behavior of $f$ on $f^{-1}(b) \subset \CC^n$. For each subset ${\J} \subset \{1,2,\ldots,n\}$ such that $\Gamma_{\infty}^{\J}(f) \supsetneq \{0\}$, let $\{\gamma_1^{\J},\gamma_2^{\J}, \ldots,\gamma_{n({\J})}^{\J}\}$ be the $(\sharp {\J}-1)$-dimensional faces of $\Gamma_{\infty}^{\J}(f)$ satisfying the condition $0 \notin \gamma_i^{\J}$. For $1 \leq i \leq n(\J)$, let $u_i^{\J} \in (\RR^{\J})^* \cap \ZZ^{\J}$ be the unique non-zero primitive vector which takes its maximum in $\Gamma_{\infty}^{\J}(f)$ exactly on $\gamma_i^{\J}$ and set 
\begin{equation}
d_i^{\J}: = \max_{v\in \Gamma_{\infty}^{\J}(f)} \langle u_i^{\J} ,v \rangle \in \ZZ_{>0}. 
\end{equation}
We call $d_i^{\J}$ the lattice distance from $\gamma_i^{\J}$ to the origin $0 \in \RR^{\J}$. For each face $\gamma_i^{\J} \prec \Gamma_{\infty}^{\J}(f)$, let $\LL(\gamma_i^{\J})$ be the smallest affine linear subspace of $\RR^n$ containing $\gamma_i^{\J}$ and $\Vol_{\ZZ}(\gamma_i^{\J}) \in \ZZ_{>0}$ the normalized $(\sharp \J -1)$-dimensional volume of $\gamma_i^{\J}$ with respect to the lattice $\ZZ^n \cap \LL(\gamma_i^{\J})$. Then by using these normalized volumes we have the following result.

\begin{thm}\label{thm:3-5}
Let $f(x) \in \CC[x_1,x_2,\ldots,x_n]$ be a polynomial on $\CC^n$. Assume that $f$ satisfies the condition ($\ast$) and is non-degenerate at infinity. Then we have
\begin{enumerate}
\item {\rm (cf. \cite{L-S})} The monodromy zeta function $\zeta_f^{\infty}(t)$ at infinity of $f$ is given by
\begin{equation}
\zeta_f^{\infty}(t)=\prod_{{\J} \colon \Gamma_{\infty}^{\J}(f) \supsetneq \{0\}}\zeta^{\infty}_{f, {\J}}(t),
\end{equation}
where for each subset $\J \subset \{1,2,\ldots,n\}$ such that $\Gamma_{\infty}^{\J}(f)\supsetneq \{0\}$ we set
\begin{equation}
\zeta^{\infty}_{f,{\J}}(t):=\prod_{i=1}^{n({\J})}(1-t^{d_i^{\J}})^{(-1)^{\sharp {\J}-1}\Vol_{\ZZ}(\gamma_i^{\J})}.
\end{equation}
\item Assume moreover that $f$ is convenient. Then for any bifurcation point $b \in B_f$ of $f$ we have
\begin{equation}
\zeta_f^b(t)=\tl{\zeta_f^b}(t).
\end{equation}
\end{enumerate}
\end{thm}

\begin{proof}
(i) Since the monodromy operators
\begin{equation}
H^j(f^{-1}(R);\CC) \simto H^j(f^{-1}(R);\CC) \hspace{10mm}(j=0,1,2,\ldots)
\end{equation}
for $R\gg 0$ are defined over $\ZZ$ and their eigenvalues are roots of unity by the monodromy theorem, $\zeta_f^{\infty}(t)$ is equal to the zeta function associated with the linear transformations
\begin{equation}
H_j(f^{-1}(R);\CC) \simto H_j(f^{-1}(R);\CC) \hspace{10mm}(j=0,1,2,\ldots)
\end{equation}
on homology groups. By the isomorphism $H_j(f^{-1}(R);\CC) \simeq H_c^{2n-2-j}(f^{-1}(R);\CC)$, we see that $\zeta_f^{\infty}(t)$ coincides with the monodromy zeta function of the constructible sheaf $Rf_!\CC_{\CC^n} \in \Dbc(\CC)$ along a sufficiently large circle $C_R\subset \CC$ ($R\gg 0$). Let us restate this result more precisely by using the nearby cycle functor. Let $j \colon \CC \longhookrightarrow \PP^1=\CC\sqcup \{\infty\}$ be the compactification and set $\F:=j_!(Rf_!\CC_{\CC^n})\in \Dbc(\PP^1)$. Take a local coordinate $h$ of $\PP^1$ in a neighborhood of $\infty \in \PP^1$ such that $\infty=\{h=0\}$. Then the monodromy zeta function of the nearby cycle $\psi_h(\F)\in \Dbc(\{\infty\})$ at the point $\infty \in \PP^1$ is nothing but $\zeta_f^{\infty}(t)$. Namely we have
\begin{equation}
\zeta_f^{\infty}(t)=\zeta_{h, \infty}(\F)(t) \in \CC(t)^*.
\end{equation}
Now let us consider $\CC^n$ as a toric variety associated with the fan $\Sigma_0$ in $\RR^n$ formed by the all faces of the first quadrant $\RR_+^n:=(\RR_{\geq 0})^n\subset \RR^n$. Let $T\simeq (\CC^*)^n$ be the open dense torus in it and for each subset ${\J}\subset \{1,2,\ldots,n\}$ denote by $T_{\J}\simeq (\CC^*)^{\sharp {\J}}$ the $T$-orbit associated with the $(n-\sharp {\J})$-dimensional face
\begin{equation}
\sigma_{\J}=\RR_+^n \cap \{u=(u_1,\ldots,u_n)\ |\ \text{$u_i=0$ for $\forall i \in {\J}$}\}
\end{equation}
of $\RR_+^n$. Then we obtain a decomposition
\begin{equation}
\CC^n =\bigsqcup_{{\J} \subset \{1,2,\ldots,n\}} T_{\J}
\end{equation}
of $\CC^n$ into $T$-orbits. Set $\F_{\J}:=j_!(Rf_!\CC_{T_{\J}})\in \Dbc(\PP^1)$. Then by a standard argument we get
\begin{equation}
\zeta_f^{\infty}(t)=\prod_{{\J} \subset \{1,2,\ldots,n\}} \zeta_{h, \infty}(\F_{\J})(t) \in \CC(t)^*.
\end{equation}
Let $f_{\J} \colon T_{\J} \simeq (\CC^*)^{\sharp {\J}} \longrightarrow \CC$ be the restriction of $f$ to $T_{\J} \subset \CC^n$. Since for ${\J} \subset \{1,2,\ldots,n\}$ such that $\Gamma_{\infty}^{\J}(f)=\{0\}$ the function $f_{\J}=f|_{T_{\J}}$ is constant, we have $\zeta_{h, \infty}(\F_{\J})(t) \equiv 1$. Hence it remains for us to show that $\zeta_{h, \infty}(\F_{\J})(t) =\zeta_{f,\J}^{\infty}(t)$ for any ${\J} \subset \{1,2,\ldots,n\}$ such that $\Gamma_{\infty}^{\J}(f)\supsetneq \{0\}$. Here we shall prove this equality only for the case where ${\J}={\J_0}:=\{1,2,\ldots,n\}$ and $T_{\J}=T = (\CC^*)^n$. The proof for the cases ${\J} \subsetneq \{1,2,\ldots,n\}$ is essentially the same. Now recall the construction of a compactification of $T \simeq (\CC^*)^n$ used in the proof of \cite[Theorem 1]{L-S}. First let $\Sigma_1$ be the dual fan of $\Gamma_{\infty}(f) \subset \RR_v^n$ in the dual vector space $\RR_u^n=(\RR_v^n)^*$ of $\RR_v^n$. Next by subdividing $\Sigma_1$ we construct a fan $\Sigma$ in $\RR_u^n$ such that the toric variety $X_{\Sigma}$ associated with it is complete and smooth. Then $T\simeq (\CC^*)^n$ is an open dense subset of $X_{\Sigma}$ and $f|_T\colon T \longrightarrow \CC$ can be extended to a meromorphic function $\tl{f}$ on $X_{\Sigma}$ which has points of indeterminacy in general. For a cone $\sigma$ in $\Sigma$ by taking a non-zero vector $u \in \relint(\sigma)$ in the relative interior $\relint(\sigma)$ of $\sigma$ we define a face $\gamma(\sigma)$ of $\Gamma_{\infty}(f)$ by
\begin{equation}
\gamma(\sigma) =\left\{ v \in \Gamma_{\infty}(f) \ \Big| \ \langle u ,v \rangle = \min_{w \in \Gamma_{\infty}(f)} \langle u,w \rangle \right\}.
\end{equation} 
By the construction of $\Sigma$, this face $\gamma(\sigma)$ does not depend on the choice of $u \in \relint(\sigma)$. Now, following \cite{L-S}, we say that a $T$-orbit in $X_{\Sigma}$ is at infinity if the corresponding cone $\sigma$ in $\Sigma$ satisfies the condition $0 \notin \gamma(\sigma)$. It is easy to see that the points of indeterminacy of $\tl{f}$ are contained the union of $T$-orbits at infinity in $X_{\Sigma}$. Let $T_1,T_2,\ldots, T_l$ be the $(n-1)$-dimensional $T$-orbits at infinity in $X_{\Sigma}$. We call their closures $D_i:=\overline{T_i}$ toric divisors at infinity in $X_{\Sigma}$. For any $i=1,2,\ldots,l$, the toric divisor $D_i$ is a smooth subvariety in $X_{\Sigma}$ and $\tl{f}$ has a pole along it. For each $i=1,2,\ldots,l$, let $\sigma_i \in \Sigma$ be the $1$-dimensional cone in $\Sigma$ which corresponds to the $T$-orbit $T_i \simeq (\CC^*)^{n-1} \subset X_{\Sigma}$. We denote the unique non-zero primitive vector in $\sigma_i \cap \ZZ^n$ by $a_i$. Then the order $m_i>0$ of the pole of $\tl{f}$ along $D_i$ is given by the formula
\begin{equation}
m_i=-\min_{v\in \Gamma_{\infty}(f)} \langle a_i,v \rangle.
\end{equation}
By the non-degeneracy of $f$ at infinity, in an open neighborhood of $\bigcup_{i=1}^lD_i$ ($=$ the union of $T$-orbits at infinity) there exists a smooth hypersurface $Z$ which satisfies the conditions:
\begin{enumerate}
\item[(a)] $\tl{f}$ has a zero of order one on $Z \setminus (\bigcup_{i=1}^lD_i)$.
\item[(b)] Any $T$-orbit at infinity in $X_{\Sigma}$ intersects $Z$ transversally.
\item[(c)] For any $T$-orbit $T_{\sigma}$ ($\sigma \in \Sigma$) at infinity in $X_{\Sigma}$, the set of points of indeterminacy of $\tl{f}$ on $T_{\sigma}$ is $T_{\sigma}\cap Z$.
\end{enumerate}
The smooth hypersurface $Z$ is uniquely determined by these conditions (a)-(c). We can show the existence of $Z$ as follows. Let $T_{\sigma}$ ($\sigma \in \Sigma$) be a $T$-orbit at infinity and $\{a_{i_1},a_{i_2},\ldots,a_{i_q},b_1,\ldots,b_r\}$ the set of the non-zero primitive vectors on the edges of $\sigma$ such that $\min_{v \in \Gamma_{\infty}(f)}\langle b_i,v\rangle =0$ for $i=1,2,\ldots,r$. Then in an open neighborhood of $T_{\sigma}$ in $X_{\Sigma}$ there exists a local coordinate system $y_1,y_2,\ldots,y_n$ such that $T_{\sigma}=\{y_1=\cdots =y_{q+r}=0, \ y_{q+r+1}, \ldots , y_n \not=0 \}$ $\simeq (\CC^*)^{n-q-r}$ and the meromorphic function $\tl{f}$ has the form
\begin{equation}
\tl{f}(y)=\dfrac{s(y)}{y_1^{m_{i_1}}y_2^{m_{i_2}} \cdots y_q^{m_{i_q}}} \hspace{10mm}(m_{i_1},\ldots, m_{i_q}>0),
\end{equation}
where $s(y)$ is a holomorphic function whose zero set is a smooth hypersurface intersecting $T_{\sigma}$ transversally. In this open set of $X_{\Sigma}$, $Z$ is explicitly defined by $Z=\{s(y)=0\}$. In order to eliminate the indeterminacy of the meromorphic function $\tl{f}$ on $X_{\Sigma}$, we first consider the blow-up $\pi_1 \colon X_{\Sigma}^{(1)} \longrightarrow X_{\Sigma}$ of $X_{\Sigma}$ along the $(n-2)$-dimensional smooth subvariety $D_1\cap Z$. Then the indeterminacy of the pull-back $\tl{f}\circ \pi_1$ of $\tl{f}$ to $X_{\Sigma}^{(1)}$ is improved. If $\tl{f}\circ \pi_1$ still has points of indeterminacy on the intersection of the exceptional divisor $E_1$ of $\pi_1$ and the proper transform $Z^{(1)}$ of $Z$, we construct the blow-up $\pi_2 \colon X_{\Sigma}^{(2)} \longrightarrow X_{\Sigma}^{(1)}$ of $X_{\Sigma}^{(1)}$ along $E_1 \cap Z^{(1)}$. By repeating this procedure $m_1$ times, we obtain a tower of blow-ups
\begin{equation}
X_{\Sigma}^{(m_1)} \underset{\pi_{m_1}}{\longrightarrow} \cdots \cdots \underset{\pi_2}{\longrightarrow} X_{\Sigma}^{(1)} \underset{\pi_1}{\longrightarrow} X_{\Sigma}.
\end{equation}
Then the pull-back of $\tl{f}$ to $X_{\Sigma}^{(m_1)}$ has no indeterminacy over $T_1$ (see the figures below). 

\bigskip
\begin{minipage}[t]{0.2\textwidth}
\begin{center}
\includegraphics[scale=1]{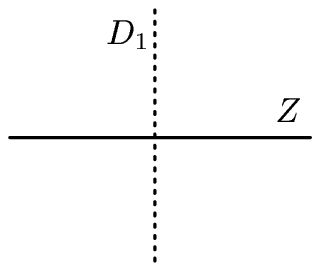}

Figure 1
\end{center}
\end{minipage}
\hspace{12mm}\begin{minipage}[t]{0.22\textwidth}
\begin{center}
\includegraphics[scale=.8]{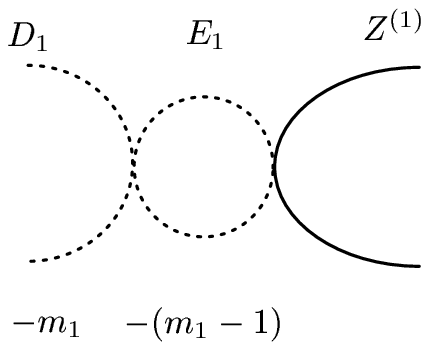}

Figure 2
\end{center}
\end{minipage}
\hspace{12mm}\begin{minipage}[t]{0.38\textwidth}
\begin{center}
\includegraphics[scale=.8]{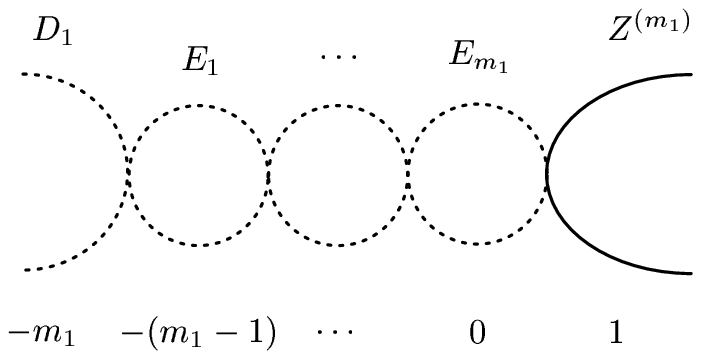}

Figure 3
\end{center}
\end{minipage}

\bigskip
Next we apply this construction to the proper transforms of $D_2$ and $Z$ in $X_{\Sigma}^{(m_1)}$. Then we obtain also a tower of blow-ups
\begin{equation}
X_{\Sigma}^{(m_1)(m_2)} \longrightarrow \cdots \cdots \longrightarrow X_{\Sigma}^{(m_1)(1)} \longrightarrow X_{\Sigma}^{(m_1)}
\end{equation}
and the indeterminacy of the pull-back of $\tl{f}$ to $X_{\Sigma}^{(m_1)(m_2)}$ is eliminated over $T_1 \sqcup T_2$. By applying the same construction to (the proper transforms of) $D_3, D_4,\ldots, D_l$, we finally obtain a birational morphism
\begin{equation}
\pi \colon \tl{X_{\Sigma}} \longrightarrow X_{\Sigma}
\end{equation}
such that $g:=\tl{f} \circ \pi$ has no point of indeterminacy on the whole $\tl{X_{\Sigma}}$ (see also \cite[p.602-604]{G-H} and \cite[Theorem 7.21]{Harris}). Then we get a commutative diagram of holomorphic maps
\begin{equation}
\xymatrix@C=20mm{ T=T_{\J_0} \ar@{^{(}->}[r]^{\iota} \ar[d]_{f|_T}& \tl{X_{\Sigma}} \ar[d]^g\\
\CC \ar@{^{(}->}[r]^j & \PP^1}
\end{equation}
and an isomorphism
\begin{eqnarray}
\F_{\J_0}
&=& j_!R(f|_T)_!\CC_T \\
&\simeq& Rg_* \iota_! \CC_T
\end{eqnarray}
in $\Dbc(\PP^1)$ ($g$ is proper). Hence in order to prove 
\begin{equation}
\zeta_{h, \infty}(\F_{\J_0})(t)=\zeta_{h, \infty}(Rg_* \iota_! \CC_T )(t)=\zeta^{\infty}_{f, \J_0}(t)
\end{equation}
we can apply Proposition \ref{prp:2-9} to the proper morphism $g \colon \tl{X_{\Sigma}} \longrightarrow \PP^1$. By calculating the monodromy zeta function of $\psi_{h \circ g}(\iota_! \CC_T)$ at each point of $(h\circ g)^{-1}(0) =g^{-1}(\infty ) \subset \tl{X_{\Sigma}}$, we can easily prove that
\begin{eqnarray} 
\zeta_{h, \infty}(\F_{\J_0})(t)
&=&\zeta_{h, \infty}(Rg_* \iota_!\CC_T)(t)\\
&=&\prod_{i=1}^l(1-t^{m_i})^{\chi(T_i\setminus Z)}.\label{eq:3-28}
\end{eqnarray}
Let us set $\delta_i:=\{v \in \Gamma_{\infty}(f)\ |\ \langle a_i, v\rangle=-m_i\}\subset \Gamma_{\infty}(f)$ for $i=1,2,\ldots,l$. Note that $\delta_i$ are (not necessarily $(n-1)$-dimensional) faces of $\Gamma_{\infty}(f)$. Then by Bernstein-Khovanskii-Kushnirenko's theorem (Theorem \ref{thm:2-10}), $\chi(T_i\setminus Z)$ is $(-1)^{n-1}$ times the normalized $(n-1)$-dimensional volume $\Vol_{\ZZ} (\delta_i)$ of $\delta_i$. Therefore we can eliminate the terms $(1-t^{m_i})^{\chi(T_i\setminus Z)}$ for $1 \leq i \leq l$ such that $\Vol_{\ZZ} (\delta_i)=0$ from \eqref{eq:3-28} and obtain the desired result
\begin{eqnarray}
\zeta_{h, \infty}(\F_{\J_0})(t)
&=& \prod_{i=1}^{n(\J_0)} \(1-t^{d_i^{\J_0}}\)^{(-1)^{\sharp \J_0-1}\Vol_{\ZZ}\(\gamma_i^{\J_0}\)}\\
&=& \zeta^{\infty}_{f,\J_0}(t).
\end{eqnarray}
This completes the proof of (i).

\noindent (ii) Recall that in the proof of (i) we identified $\CC^n$ with the toric variety associated with the fan $\Sigma_0$. Since $f$ is convenient, $\Sigma_0$ is an open subfan of $\Sigma_1$. Then, without subdividing the cones in $\Sigma_0\subset \Sigma_1$ we can construct a subdivision $\Sigma$ of $\Sigma_1$ such that the associated toric variety $X_{\Sigma}$ is complete and smooth. Namely $\CC^n$ is an affine open subset of $X_{\Sigma}$. Let $\iota^{\prime} \colon \CC^n \longhookrightarrow X_{\Sigma}$ be the inclusion. Then by eliminating the points of indeterminacy of the meromorphic extension $\tl{f}$ of $f$ to $X_{\Sigma}$ along toric divisors at infinity in $X_{\Sigma}$ as in (i) we obtain a commutative diagram of holomorphic maps
\begin{equation}
\xymatrix{
\CC^n \ar@{^{(}->}[r]^{\iota^{\prime}} \ar[d]_f&\tl{X_{\Sigma}} \ar[d]^{g^{\prime}}\\
\CC \ar@{^{(}->}[r]^j & \PP^1}
\end{equation}
and an isomorphism
\begin{eqnarray}
\F&=& j_! Rf_!\CC_{\CC^n} \\
&\simeq& Rg^{\prime}_* \iota^{\prime}_!\CC_{\CC^n}.
\end{eqnarray}
Now let $h^{\prime}$ be a local coordinate of $\PP^1$ on an open neighborhood of the bifurcation point $b\in \PP^1$ such that $b=\{h^{\prime}=0\}$. Then we have
\begin{eqnarray}
\zeta_f^b(t)
&=& \zeta_{h^{\prime}, b}(\F)(t)\\
&=& \zeta_{h^{\prime}, b}(Rg^{\prime}_* \iota^{\prime}_!\CC_{\CC^n})(t).
\end{eqnarray}
Therefore, as in the proof of (i), by applying Proposition \ref{prp:2-9} to the proper morphism $g^{\prime} \colon \tl{X_{\Sigma}} \longrightarrow \PP^1$ we can easily prove the desired result
\begin{equation}
\zeta_f^b(t)=\tl{\zeta_f^b}(t).
\end{equation}
This completes the proof.\qed
\end{proof}

\begin{cor}
Let $f(x) \in \CC[x_1,x_2,\ldots,x_n]$ be a polynomial on $\CC^n$. Assume that $f$ is convenient and non-degenerate at infinity. Let $b \in B_f$ be a bifurcation point of $f$ such that the complex hypersurface $f^{-1}(b) \subset \CC^n$ has only isolated singular points $\{p_1,p_2,\ldots,p_k\} \subset f^{-1}(b)$. Denote by $\zeta_i(t)$ the local monodromy zeta function $\zeta_{f-b, p_i}(t)$ of the complex hypersurface $f^{-1}(b)=\{ f-b=0 \}$ at $p_i$. Then we have
\begin{equation}
\zeta_f^b(t)=\tl{\zeta_f^b}(t)=(1-t)^{\chi(f^{-1}(b)\setminus\{p_1,p_2,\ldots,p_k\})} \prod_{i=1}^k \zeta_i(t).
\end{equation}
\end{cor}

\begin{rem}
For $k\in \NN$, we denote the Lefschetz number of $(\Phi_f^{\infty})^k$ by $L^k(\Phi_f^{\infty})$. Since we have
\begin{equation}
\zeta_f^{\infty}(t) =\exp \( -\sum_{k=1}^{\infty} \dfrac{L^k(\Phi_f^{\infty})}{k}t^k\),
\end{equation}
we can calculate the Lefschetz number $L^k(\Phi_f^{\infty})$ as follows. Namely, in the situation of Theorem \ref{thm:3-5} we have
\begin{eqnarray}
L^k(\Phi_f^{\infty})
&=&\dfrac{-1}{(k-1)!}\( \dfrac{d}{dt}\)^k\( \log{\(\zeta_f^{\infty}(t)\)}\)\Big|_{t=0}\\
&=&\sum_{\J \colon \Gamma_{\infty}^{\J}(f) \supsetneq \{0\}} \sum_{i=1}^{n(\J)}\dfrac{(-1)^{\sharp \J}\Vol_{\ZZ}(\gamma_i^{\J})}{(k-1)!}\( \dfrac{d}{dt}\)^k\(\log{(1-t^{d_i^{\J}})}\)\Big|_{t=0}\\
&=&\sum_{\J \colon \Gamma_{\infty}^{\J}(f) \supsetneq \{0\}}\sum_{\begin{subarray}{c}0\leq i \leq n(\J)\\ d_i^{\J}|k\end{subarray}}(-1)^{\sharp \J-1}\Vol_{\ZZ}(\gamma_i^{\J})d_i^{\J}.
\end{eqnarray}
By this formula we observe that $L^k(\Phi_f^{\infty})$ does not always vanish (compare this result with the one obtained by A'Campo in \cite{A'Campo}). Indeed in the case where $f(x_1,x_2)=x_1(x_1^2x_2^2-1)$ we have $L^k(\Phi_f^{\infty})=-1$ for any $k\in \NN$. Moreover, if we set $Z_f(r):=\dfrac{d}{dr}\log{(\zeta_f^{\infty}(e^r))}$, then we obtain also a functional equation
\begin{eqnarray}
Z_f(r)+Z_f(-r)
&=&\sum_{\J \colon \Gamma_{\infty}^{\J}(f) \supsetneq \{0\}} \sum_{i=1}^{n(\J)}(-1)^{\sharp \J-1}\Vol_{\ZZ}(\gamma_i^{\J})d_i^{\J}\\
&=&\deg_t \zeta_f^{\infty}(t)=\chi(f^{-1}(R)).
\end{eqnarray}
We can obtain similar formulas also for other topological zeta functions.
\end{rem}

\section{Global monodromy along central fibers}\label{sec:4}

From now on, we shall generalize Theorem \ref{thm:3-5} (ii) to non-convenient polynomials. In Theorem \ref{thm:4-5} and \ref{thm:4-9} below, we obtain unexpected results that the constant term $a=a_0\in \CC$ of a non-convenient polynomial $f(x)=\sum_{v \in \ZZ_+^n} a_v x^v$ ($a_v\in \CC$) on $\CC^n$ often becomes a bifurcation point of $f$. In what follows, we assume always that the dimension of the Newton polygon at infinity $\Gamma_{\infty}(f)$ of $f(x)\in \CC[x_1,\ldots,x_n]$ is $n$. We consider the $n$-dimensional cone $\Cone_{\infty}(f)=\RR_{\geq 0}\Gamma_{\infty}(f)$ generated by $\Gamma_{\infty}(f)$ in $\RR^n$.

\begin{dfn}\label{dfn:4-2}
For a polynomial $f(x)=\sum_{v \in \ZZ_+^n}a_v x^v$ ($a_v \in \CC$) on $\CC^n$ with the constant term $a:=a_0 \in \CC$, we set
\begin{equation}
\Gamma_{\bif}(f):=NP(f-a)+\Cone_{\infty}(f) \subset \RR^n.
\end{equation}
We call $\Gamma_{\bif}(f)$ the bifurcation Newton polygon of $f$ at $a$.
\end{dfn}

If $n=2$, the following definition is useful. Let $\RR_u^n=(\RR_v^n)^*$ be the dual vector space of $\RR_v^n$ and consider the first quadrant $\RR_+^n=\RR_{\geq 0}^n \subset \RR_u^n$ in it and its interior $\Int(\RR_+^n)=\RR_{> 0}^n$.

\begin{dfn}\label{dfn:4-4}
Let $f(x)=\sum_{v \in \ZZ_+^n}a_vx^v$ ($a_v \in \CC$) be a polynomial on $\CC^n$ and $c \in \CC$. Then we say $f$ is non-degenerate (resp. strongly non-degenerate) along the fiber $f^{-1}(c)$ if for any $u \in \RR_u^n \setminus \RR_+^n$ (resp. for any non-zero $u \in \RR_u^n \setminus \Int(\RR_+^n)$) the complex hypersurface
\begin{equation}
\{x=(x_1,\ldots,x_n)\in (\CC^*)^n \ |\ (f-c)^u(x)=0\}
\end{equation}
in $(\CC^*)^n$ is smooth and reduced.
\end{dfn}

Note that in the above definition we do not assume that the dimension of $NP(f-c)$ is $n$.

\begin{thm}\label{thm:4-5}
Let $f(x)=\sum_{v \in \ZZ_+^2}a_vx^v$ ($a_v\in \CC$) be a polynomial on $\CC^2$ with the constant term $a:=a_0 \in \CC$.
\begin{enumerate}
\item Assume that $f$ is non-degenerate along the fiber $f^{-1}(a) $. Then we have
\begin{equation}\label{eq:4-3}
\zeta_f^a (t)=\prod_{j=1}^l (1-t^{d_j})^{-\Vol_{\ZZ}(\gamma_j)}\times \tl{\zeta_f^a }(t),
\end{equation}
where $\gamma_1,\gamma_2,\ldots,\gamma_l$ are the $1$-dimensional compact faces of $\Gamma_{\bif}(f)$ whose inner conormal vectors are contained in $\RR_u^2\setminus \RR_+^2$, $d_j\in \ZZ_{>0}$ is the lattice distance from $\gamma_j$ to the origin $0 \in \RR_v^2$ and $\Vol_{\ZZ}(\gamma_j)$ is the normalized volume of $\gamma_j$ with respect to the lattice $\ZZ^2 \cap \LL(\gamma_j)$.
\item For a complex number $c\neq a$, assume that $f$ is non-degenerate along the fiber $f^{-1}(c)$. Then we have
\begin{equation}
\zeta_f^c(t)= \tl{\zeta_f^c}(t).
\end{equation}
\end{enumerate}
\end{thm}

\begin{minipage}[t]{0.4\textwidth}
\begin{center}
\includegraphics[scale=.7]{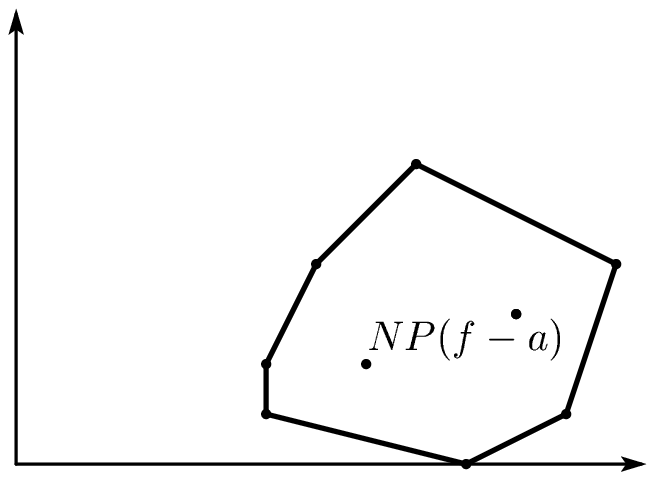}

Figure 4
\end{center}
\end{minipage}
\hspace{20mm}\begin{minipage}[t]{0.4\textwidth}

\begin{center}
\includegraphics[scale=.7]{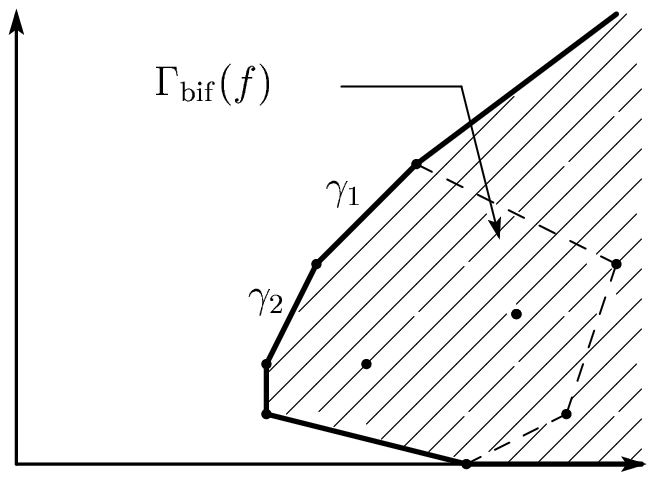}

Figure 5
\end{center}
\end{minipage}

\begin{proof}
Since the proof of (ii) is similar, we prove only (i). As in the proof of Theorem \ref{thm:3-5}, let $\Sigma_0$ be the fan in $\RR_u^2$ formed by faces of the first quadrant $\RR_+^2$. We also denote by $\Sigma_1$ the dual fan of the $2$-dimensional polytope $\Gamma_{\infty}(f) \subset \RR_v^2$. For $j=1,2,\ldots,l$, let $u_j \in \RR_u^2 \setminus \RR_+^2$ be the primitive inner conormal vector of $\gamma_j$ in $\ZZ^2\subset \RR_u^2$ and consider the rays (i.e. $1$-dimensional cone) $r_j =\RR_{\geq 0}u_j$ in $\RR_u^2$ generated by $u_j$. By subdividing $\Sigma_1$ by $r_1,r_2,\ldots, r_l$, we obtain a fan $\Sigma_2$ in $\RR_u^2$. Then there exists a subdivision $\Sigma$ of $\Sigma_2$ which contains the unique $2$-dimensional cone in $\Sigma_0$ such that the associated toric variety $X_{\Sigma}$ is smooth. In this situation, $\CC^2$ is an affine open subset of the smooth toric surface $X_{\Sigma}$. Let $T \simeq (\CC^*)^2$ be the open dense torus in $X_{\Sigma}$ which acts $X_{\Sigma}$ itself. Then by the non-degeneracy of $f$ along the fiber $f^{-1}(a)$, the closure $\overline{C}$ of the complex curve $C=\{ x\in \CC^2 \ | \ f(x)-a=0 \} \subset \CC^2$ in $X_{\Sigma}$ intersects each $T$-orbit in $X_{\Sigma} \setminus \CC^2$ transversally. Now, as in the proof of Theorem \ref{thm:3-5}, by eliminating the points of indeterminacy of the meromorphic extension $\tl{f}$ of $f$ to $X_{\Sigma}$, we construct a complete surface $\tl{X_{\Sigma}}$ and the following commutative diagram of holomorphic maps:
\begin{equation}
\xymatrix{ \CC^2 \ar@{^{(}->}[r]^{\iota} \ar[d]_{f}& \tl{X_{\Sigma}} \ar[d]^g\\
\CC \ar@{^{(}->}[r]^j & \PP^1.}
\end{equation}
Let $h$ be a local coordinate of $\PP^1$ in a neighborhood of $a \in \PP^1$ such that $a =\{ h=0\}$. Then by using the properness of $g$ we obtain
\begin{equation}
\zeta_f^a(t)=\zeta_{h, a}(j_!Rf_! \CC_{\CC^2})(t)=\zeta_{h, a}(Rg_* \iota_! \CC_{\CC^2})(t).
\end{equation}
By Proposition \ref{prp:2-9}, for the calculation of $\zeta_{h, a}(Rg_* \iota_! \CC_{\CC^2})(t)$ it suffices to calculate the monodromy zeta function of $\psi_{h \circ g}(\iota_! \CC_{\CC^2})$ at each point of $(h \circ g)^{-1}(0)=g^{-1}(a) \subset \tl{X_{\Sigma}}$. If we denote by $T_j$ ($\simeq \CC^*$) the $1$-dimensional $T$-orbit in $X_{\Sigma}\setminus \CC^2$ which corresponds to the ray $r_j$, we can easily see that
\begin{equation}
\zeta_{h, a}(Rg_* \iota_! \CC_{\CC^2})(t) =\prod_{j=1}^l (1-t^{d_j})^{\chi(T_j \setminus \overline{C})} \times \tl{\zeta_f^a }(t). 
\end{equation}
Since $\chi(T_j \setminus \overline{C})=- \Vol_{\ZZ}(\gamma_j)$ by Bernstein-Khovanskii-Kushnirenko's theorem (Theorem \ref{thm:2-10}), we obtain the desired formula \eqref{eq:4-3}. \qed
\end{proof}

\begin{exa}[\cite{Broughton}]
Let $f(x_1, x_2)$ be a polynomial on $\CC^2$ defined by
\begin{equation}
f(x_1, x_2)=x_1-x_1^2x_2=x_1(1-x_1x_2). 
\end{equation}
In this case, the constant term $a$ of $f$ is $0$ and by using the homeomorphism $f^{-1}(0) \simeq \CC \sqcup \CC^*$ we have $\tl{\zeta_f^0 }(t)=(1-t)$. Then by Theorem \ref{thm:4-5} we obtain 
\begin{equation}
\zeta_f^0 (t)=(1-t)^{-1} \times \tl{\zeta_f^0 }(t)=1
\end{equation}
(see also \cite[Example 6.1]{S-T-2}). Since $0 \in \CC$ is the unique bifurcation point of the polynomial map $f\colon \CC^2 \longrightarrow \CC$ in this case, we have $\zeta_f^0 (t)=\zeta_f^{\infty}(t)$ and the equality $\zeta_f^{0}(t)=1$ can be deduced also from Theorem \ref{thm:3-5} (i).
\end{exa}

Since for the constant term $a=a_0 \in \CC$ of a non-convenient polynomial $f(x)=\sum_{v \in \ZZ_+^2}a_vx^v$ on $\CC^2$ the global fiber $f^{-1}(a) \subset \CC^2$ may have multiplicities, we regard $f^{-1}(a)$ as a reduced divisor and let $p_1,p_2,\ldots, p_k$ be the isolated singular points of $f^{-1}(a)$. We set $\mu_i := 1-\deg(\zeta_{f-a,p_i}(t)) \in \ZZ$. Then by using the fact that the degree of $\zeta_f^a(t)$ in $t$ is the Euler characteristic of the general fiber $f^{-1}(c)$ ($c\in \CC \setminus B_f$) we can easily prove the following formula for the jumping number
\begin{equation}
\chi(f^{-1}(a))-\chi(f^{-1}(c))
\end{equation}
of the Euler characteristic of $f^{-1}(a)$ from that of the general fiber $f^{-1}(c)$.

\begin{cor}\label{cor:4-6}
Let $f(x)=\sum_{v \in \ZZ_+^2}a_vx^v$ ($a_v\in \CC$) be a polynomial on $\CC^2$ with the constant term $a=a_0 \in \CC$. Assume that $f$ is non-convenient and strongly non-degenerate along the fiber $f^{-1}(a)$. We define two non-negative integers $m_1, m_2 \in\ZZ_+$ by
\begin{gather}
m_1:=\min \{ v_1 \ | \ v=(v_1, v_2) \in NP(f-a) \}, \\
m_2 :=\min \{ v_2 \ | \ v=(v_1, v_2) \in NP(f-a) \}.
\end{gather}
Denote by $L_x \in \ZZ_+$ (resp. $L_y \in \ZZ_+$) the length (normalized $1$-dimensional volume) of the segment $NP(f-a) \cap \{ v_2=m_2 \}$ (resp. $NP(f-a) \cap \{ v_1=m_1 \}$). Then for any $c\in \CC \setminus B_f$ we have
\begin{equation}
\chi(f^{-1}(a))-\chi(f^{-1}(c))=\sum_{j=1}^l d_j\cdot \Vol_{\ZZ}(\gamma_j) +\sum_{i=1}^k \mu_i +K,
\end{equation}
where the last term $K$ is defined by
\begin{equation}
K:=\begin{cases}
(m_2 -1)L_x+(m_1 -1)L_y &(m_1, \ m_2 >0),\\
(m_1 -1)L_y &(m_1 >0, \ m_2 =0, \ (m_1, m_2) \notin \supp (f)),\\
(m_1 -1)(L_y-1) &(m_1 >0, \ m_2 =0, \ (m_1, m_2) \in \supp (f)),\\
(m_2 -1)L_x &(m_1 =0, \ m_2 >0, \ (m_1, m_2) \notin \supp (f)),\\
(m_2 -1)(L_x-1) &(m_1 =0, \ m_2 >0, \ (m_1, m_2) \in \supp (f)).
\end{cases}
\end{equation}
\end{cor}

\begin{proof}
By the strong non-degeneracy of $f$ along the fiber $f^{-1}(a)$, the closure $\overline{C}$ of the complex curve $C=\{ x \in (\CC^*)^2 \ | \ f(x)-a=0\} \subset (\CC^*)^2$ in $\CC^2$ intersects $1$-dimensional tori $\CC^* \times \{ 0\}$, $\{ 0\} \times \CC^* \subset \CC^2$ transversally. Moreover, by Bernstein-Khovanskii-Kushnirenko's theorem (Theorem \ref{thm:2-10}), we obtain 
\begin{equation}
\sharp \{ \overline{C} \cap (\CC^* \times \{ 0\})\} =L_x, \qquad \sharp \{ \overline{C} \cap (\{ 0\} \times \CC^*)\} =L_y.
\end{equation} 
Then the result follows easily from Theorem \ref{thm:4-5}. \qed
\end{proof}

To extend Theorem \ref{thm:4-5} to the case $n\geq 3$, we need the following definition.

\begin{dfn}\label{dfn:4-7}
Let $f(x)=\sum_{v \in \ZZ_+^n} a_vx^v$ ($a_v \in \CC$) be a polynomial on $\CC^n$ with the constant term $a:=a_0 \in \CC$. For ${\J} \subset \{1,2,\ldots,n\}$ such that $\Gamma_{\infty}^{\J}(f)\supsetneq \{0\}$ ($\Longleftrightarrow$ $NP(f-a) \cap \RR^{\J} \neq \emptyset$), we set
\begin{equation}
\Gamma_{\bif}^{\J}(f):=\Gamma_{\bif}(f) \cap\RR^{\J} \subset \RR^{\J}.
\end{equation}
\end{dfn}

The following definition is just a slight a modification of the standard one.

\begin{dfn}\label{dfn:4-8}
Let $f(x)=\sum_{v \in \ZZ_+^n} a_v x^v $ ($a_v \in \CC$) be a polynomial on $\CC^n$ and $c \in \CC$. Then we say that $f$ is strictly non-degenerate along the fiber $f^{-1}(c)$ if for any non-zero $u \in (\RR^n)^*$ the complex hypersurface
\begin{equation}
\{x=(x_1,x_2,\ldots,x_n)\in (\CC^*)^n \ |\ (f-c)^{u}(x)=0\}
\end{equation}
in $(\CC^*)^n$ is smooth and reduced.
\end{dfn}

For each ${\J} \subset \{1,2,\ldots,n\}$, consider the algebraic torus
\begin{equation}
T_{\J}:=\{x=(x_1,\ldots,x_n) \in \CC^n \ |\ x_i=0\ (i \notin {\J}), \ x_i\neq 0\ (i \in {\J})\}\simeq (\CC^*)^{\sharp {\J}}
\end{equation}
and the standard decomposition $\CC^n =\bigsqcup_{{\J} \subset \{1,2,\ldots,n\}} T_{\J}$ of $\CC^n$. Let $f_{\J} \colon T_{\J} \simeq (\CC^*)^{\sharp {\J}} \longrightarrow \CC$ be the restriction of $f$ to $T_{\J} \subset \CC^n$. Then for $c \in \CC$ and a subset ${\J} \subset \{1,2,\ldots,n\}$ such that $\Gamma_{\infty}^{\J}(f) \supsetneq \{0\}$ ($\Longleftrightarrow$ $f_{\J}-c$ is not constant), by taking the Euler integral of the local monodromy zeta function $\zeta_{f_{\J}-c}(t)$ over $\{ f_{\J}(x)-c=0 \} \subset T_{\J}$ we obtain a rational function $\tl{\zeta_{f,{\J}}^c}(t) \in \CC(t)^*$.

\begin{thm}\label{thm:4-9}
Let $f(x)=\sum_{v \in \ZZ_+^n}a_vx^v$ ($a_v \in \CC$) be a polynomial on $\CC^n$ with the constant term $a:=a_0 \in \CC$.
\begin{enumerate}
\item Assume that $f$ is strictly non-degenerate along the fiber $f^{-1}(a)$. Then the monodromy zeta function $\zeta_f^a(t)$ is given by
\begin{equation}
\zeta_f^a(t)=\prod_{{\J} \colon \Gamma_{\infty}^{\J}(f) \supsetneq \{0\}} \zeta_{f,{\J}}^a(t),
\end{equation}
where for each ${\J} \subset \{1,2,\ldots,n\}$ such that $\Gamma_{\infty}^{\J}(f) \supsetneq \{0\}$ we set
\begin{equation}
\zeta_{f,{\J}}^a(t):=\prod_{j=1}^{\nu ({\J})} (1-t^{e_j^{\J}})^{(-1)^{\sharp {\J}-1} \Vol_{\ZZ}(\delta_j^{\J})}\times \tl{\zeta_{f,{\J}}^a}(t).
\end{equation}
Here $\{\delta_1^{\J},\delta_2^{\J},\ldots,\delta_{\nu ({\J})}^{\J}\}$ is the set of $(\sharp {\J}-1)$-dimensional compact faces of $\Gamma_{\bif}^{\J}(f)$, $e_j^{\J} \in \ZZ_{>0}$ is the lattice distance from $\delta_j^{\J}$ to the origin $0\in \RR^{\J}$ and $\Vol_{\ZZ}(\delta_j^{\J})$ is the normalized $(\sharp {\J}-1)$-dimensional volume of $\delta_j^{\J}$ with respect to the lattice $\ZZ^n \cap \LL(\delta_j^{\J})$.
\item For a complex number $c \neq a$, assume that $f$ is strictly non-degenerate along the fiber $f^{-1}(c)$. Then we have
\begin{equation}
\zeta_{f}^c(t)=\tl{\zeta_f^c}(t).
\end{equation}
\end{enumerate}
\end{thm}

\begin{proof}
(i) First note that for ${\J} \subset \{1,2,\ldots,n\}$ such that $\Gamma_{\infty}^{\J}(f)=\{0\}$ the function $f_{\J}-a$ on $T_{\J} \simeq (\CC^*)^{\sharp {\J}}$ is identically zero. Let $h$ be a local coordinate of $\CC$ in a neighborhood of $a \in \CC$ such that $a=\{ h=0\}$. Then as in the proof of Theorem \ref{thm:3-5} (i), it is enough to prove that for any ${\J} \subset \{1,2,\ldots,n\}$ such that $\Gamma_{\infty}^{\J}(f) \supsetneq \{0\}$ we have 
\begin{equation}\label{eq:4-22} 
\zeta_{h, a}(R(f_{\J})_! \CC_{T_{\J}})(t)=\zeta_{f,{\J}}^a(t). 
\end{equation}
We prove this equality only for the case where ${\J}=\J_0:=\{1,2,\ldots,n\}$. Let $T:=T_{\J_0} =(\CC^*)^n$ be the open dense torus in $\CC^n$. Let $\Sigma_1$ be the dual subdivision of $NP(f-a)$ in $\RR_u^n=(\RR_v^n)^*$. Note that $\Sigma_1$ is not necessarily a fan in $\RR_u^n$ since here we do not assume that the dimension of $NP(f-a)$ is $n$. By dividing $\Sigma_1$ into strongly convex rational cones and applying some more subdivisions if necessary, we obtain a fan $\Sigma$ in $\RR_u^n$ such that the associated toric variety $X_{\Sigma}$ is smooth and complete. Then $T=(\CC^*)^n$ is an open dense subset of $X_{\Sigma}$. By eliminating the points of indeterminacy of the meromorphic extension $\tl{f}$ of $f$ to $X_{\Sigma}$ as in the proof of Theorem \ref{thm:3-5} (i), we obtain a complete variety $\tl{X_{\Sigma}}$ and the following commutative diagram of holomorphic maps. 
\begin{equation}
\xymatrix{
T \ar@{^{(}->}[r]^{\iota} \ar[d]_{f_{\J_0}=f|_T}& \tl{X_{\Sigma}} \ar[d]^g\\ \CC \ar@{^{(}->}[r]^j& \PP^1.}
\end{equation}
By Proposition \ref{prp:2-9}, for the calculation of $\zeta_{h,a}(R(f_{\J_0})_!\CC_{T})(t)=\zeta_{h,a}(Rg_* \iota_! \CC_{T})(t)$, it suffices to calculate the monodromy zeta function of $\psi_{h \circ g}(\iota_! \CC_{T})(t)$ at each point of $(h \circ g)^{-1}(0)=g^{-1}(a) \subset \tl{X_{\Sigma}}$. Then, as the proof of Theorem \ref{thm:3-5} (i), we can easily prove the equality \eqref{eq:4-22}.

\noindent (ii) In the same way as the proof of (i), for $c \neq a$ we can prove that
\begin{equation}
\zeta_{f}^c(t)=\prod_{{\J} \colon \Gamma_{\infty}^{\J}(f) \supsetneq \{0\}} \tl{\zeta_{f,{\J}}^c} (t).
\end{equation}
By the strict non-degeneracy of $f$ along the fiber $f^{-1}(c)$, the hypersurface $\{ f(x)-c=0 \}$ in $\CC^n$ intersects $T_{\J}$ transversally for any ${\J} \subset \{ 1,2, \ldots , n\}$ such that $\Gamma_{\infty}^{\J}(f) \supsetneq \{0\}$. Then for such $\J$ the local monodromy zeta functions $\zeta_{f_{\J}-c, x}(t)$ and $\zeta_{f-c,x}(t)$ coincide each other at any point $x$ of $\{ f_{\J}(x)-c=0 \} \subset T_{\J}$. Moreover, for ${\J} \subset \{ 1,2, \ldots , n\}$ such that $\Gamma_{\infty}^{\J}(f) = \{0\}$ the function $f_{\J}-c=f|_{T_{\J}}-c$ on $T_{\J}$ is identically equal to the ``non-zero" complex number $a-c$. This implies that the hypersurface $\{ x\in \CC^n \ | \ f(x)-c=0 \} \subset \CC^n$ does not intersect such $T_{\J}$. Summarizing these arguments we obtain
\begin{eqnarray}
\tl{\zeta_{f}^c} (t)
&=&\prod_{{\J} \colon \Gamma_{\infty}^{\J}(f) \supsetneq \{0\}} \tl{\zeta_{f,{\J}}^c} (t)\\
&=& \zeta_{f}^c(t).
\end{eqnarray}
This completes the proof.\qed
\end{proof}

Now let $f(x)=\sum_{v \in \ZZ_+^n}a_vx^v$ ($a_v \in \CC$) be a polynomial on $\CC^n$ with the constant term $a:=a_0 \in \CC$ and $c \in \CC$. Assume that $f$ is strictly non-degenerate along the fiber $f^{-1}(c)$. Then by the proof of Theorem \ref{thm:4-9} (ii), for any ${\J} \subsetneq \J_0:= \{1,2,\ldots,n\}$ such that $\Gamma_{\infty}^{\J}(f) \supsetneq \{0\}$ the complex hypersurface $f_{\J}^{-1}(c)$ in $T_{\J}\simeq (\CC^*)^{\sharp {\J}}$ is smooth. 

\begin{cor}\label{cor:4-10}
In the situation as above, assume moreover that the complex hypersurface $f_{\J_0}^{-1}(c)$ in $T:=T_{\J_0}\simeq (\CC^*)^n$ has only isolated singular points $p_1,p_2,\ldots, p_{k} \in f_{\J_0}^{-1}(c)$. For $1 \leq i \leq k$ denote by $\zeta_i(t)$ (resp. $\mu_i \geq 1$) the local monodromy zeta function (resp. the local Milnor number) of $f_{\J_0}^{-1}(c)$ at $p_i$. Then we have
\begin{enumerate}
\item If $c=a$, the monodromy zeta function $\zeta_f^a(t)$ is given by
\begin{equation}
\zeta_f^a(t)=\prod_{{\J} \colon \Gamma_{\infty}^{\J}(f)\supsetneq \{0\}} \zeta_{f,{\J}}^a(t),
\end{equation}
where we set
\begin{equation}
\zeta_{f,{\J_0}}^a(t)
:=\prod_{j=1}^{\nu({\J_0})} (1-t^{e_j^{\J_0}})^{(-1)^{\sharp {\J_0}-1} \Vol_{\ZZ}(\delta_j^{\J_0})}\times (1-t)^{\chi(f_{\J_0}^{-1}(a) \setminus \{ p_1,p_2, \ldots , p_k\})} \times \prod_{i=1}^{k} \zeta_i(t)
\end{equation}
and for each ${\J} \subsetneq \J_0= \{1,2,\ldots,n\}$ such that $\Gamma_{\infty}^{\J}(f) \supsetneq \{0\}$ we set
\begin{equation}
\zeta_{f,{\J}}^a(t):=\prod_{j=1}^{ \nu ({\J})} (1-t^{e_j^{\J}})^{(-1)^{\sharp {\J}-1} \Vol_{\ZZ}(\delta_j^{\J})}\times (1-t)^{\chi(f_{\J}^{-1}(a))}. 
\end{equation}
Here $\{\delta_1^{\J},\delta_2^{\J}, \ldots, \delta_{\nu({\J})}^{\J}\}$ and $e_j^{\J} \in \ZZ_{>0}$ etc. are the ones used in Theorem \ref{thm:4-9} (i).
\item If $c \neq a$, the monodromy zeta function $\zeta_f^c(t)$ is given by
\begin{equation}
\zeta_f^c(t)
=(1-t)^{\chi(f_{\J_0}^{-1}(c)\setminus \{p_1,p_2,\ldots,p_k\})}\times \prod_{i=1}^{k} \zeta_i(t) \times \prod_{\J \subsetneq \J_0 \colon \Gamma_{\infty}^{\J}(f) \supsetneq \{0\}} (1-t)^{\chi(f_{\J}^{-1}(c))}. 
\end{equation}
\item If $c=a$, for any $c^{\prime} \in \CC \setminus B_f$ we have
\begin{equation}
\chi(f^{-1}(a))-\chi(f^{-1}(c^{\prime}))
= (-1)^n \sum_{i=1}^{k} \mu_i + \sum_{{\J} \colon \Gamma_{\infty}^{\J}(f) \supsetneq\{0\}} (-1)^{\sharp {\J}} \left\{\sum_{j=1}^{\nu ({\J})} e_j^{\J}\cdot \Vol_{\ZZ}(\delta_j^{\J})\right\} +1.
\end{equation}
\end{enumerate}
\end{cor}

In many cases we can easily rewrite the result of Theorem \ref{thm:4-9} (i) in terms of the finite part $\tl{\zeta_f^a}(t)$ of $\zeta_f^a(t)$. For this purpose, we introduce the following definition.

\begin{dfn}
Let $h(x)=\sum_{v \in \ZZ_+^n} a_vx^v \in \CC[x_1,x_2,\ldots, x_n]$ ($a_v \in \CC$) be a polynomial on $\CC^n$. We associate to it a polynomial $\tl{h}(x) \in \CC[x_1,x_2,\ldots, x_n]$ defined by
\begin{equation}
\tl{h}(x):=\dfrac{h(x)}{x_1^{m_1}x_2^{m_2}\cdots x_n^{m_n}},
\end{equation}
where for $1 \leq i \leq n$ we set
\begin{equation}
m_i:=\min \{v_i \ | \ v=(v_1,v_2, \ldots , v_n) \in NP(h) \}.
\end{equation}
\begin{enumerate}
\item We say that $h$ is semi-convenient if $\tl{h}$ is convenient.
\item We say that $h$ is quasi-convenient if for any face $\sigma \prec \RR_+^n$ of $\RR_+^n$ such that $\d \sigma \geq 1$ we have $NP(\tl{h}) \cap \sigma \neq \emptyset$.
\end{enumerate}
\end{dfn}

Note that by definition semi-convenient polynomials are quasi-convenient. Moreover any polynomial on $\CC^2$ is quasi-convenient.

\begin{thm}\label{thm:4-12}
Let $f(x)=\sum_{v \in \ZZ_+^n}a_vx^v$ ($a_v \in \CC$) be a polynomial on $\CC^n$ with the constant term $a:=a_0 \in \CC$. Assume that $f-a$ is quasi-convenient and $f$ is strictly non-degenerate along the fiber $f^{-1}(a)$. Then we have
\begin{equation}
\zeta_f^a(t)=\prod_{{\J} \colon \Gamma_{\infty}^{\J}(f) \supsetneq \{0\}} \left\{ \prod_{j=1}^{l({\J})} (1-t^{d_j^{\J}})^{(-1)^{\sharp {\J}-1} \Vol_{\ZZ}(\gamma_j^{\J})}\right\} \times \tl{\zeta_f^a}(t),
\end{equation}
where for each ${\J} \subset \{1,2,\ldots,n\}$ such that $\Gamma_{\infty}^{\J}(f) \supsetneq \{0\}$, $\{\gamma_1^{\J},\gamma_2^{\J},\ldots,\gamma_{l({\J})}^{\J}\}$ is the set of the $(\sharp {\J}-1)$-dimensional compact faces of $\Gamma_{\bif}^{\J}(f)$ whose inner conormal vectors are contained in $\RR^{\J} \setminus \RR_+^{\J}$ and $d_j^{\J} \in \ZZ_{>0}$ is the lattice distance from $\gamma_j^{\J}$ to the origin $0\in \RR^{\J}$.
\end{thm}

\begin{proof}
First, by Theorem \ref{thm:4-9} we have 
\begin{equation}\label{TK1} 
\zeta_f^a(t)=\prod_{{\J} \colon \Gamma_{\infty}^{\J}(f) \supsetneq \{0\}} \left\{\prod_{j=1}^{\nu ({\J})} (1-t^{e_j^{\J}})^{(-1)^{\sharp {\J}-1} \Vol_{\ZZ}(\delta_j^{\J})}\times \tl{\zeta_{f, {\J} }^a}(t) \right\}, 
\end{equation}
where $\{\delta_1^{\J},\ldots,\delta_{\nu ({\J})}^{\J}\}$ is the set of $(\sharp {\J}-1)$-dimensional compact faces of $\Gamma_{\bif}^{\J}(f)$ and $e_j^{\J} \in \ZZ_{>0}$ is the lattice distance from $\delta_j^{\J}$ to the origin $0\in \RR^{\J}$. As we saw in the proof of Theorem \ref{thm:4-9} (ii), the strict non-degeneracy of $f$ along the fiber $f^{-1}(a)$ implies that
\begin{equation}
\int_{T_{\J}} \zeta_{f-a}(t)=\tl{\zeta_{f,{\J}}^a}(t)
\end{equation}
for any $\J \subset \{1,2,\ldots,n\}$ such that $\Gamma_{\infty}^{\J}(f) \supsetneq \{0\}$. For such $\J \subset \{1,2,\ldots,n\}$, we shall define three subsets $I_1^{\J}, I_2^{\J}, I_3^{\J}$ of $\{1,2, \ldots,\nu(\J)\}$ as follows. For $1 \leq j \leq \nu (\J)$, let $u_j^{\J} \in (\RR^{\J})^* \cap \ZZ^{\J}$ be the (unique) non-zero primitive vector which takes its minimum in $\Gamma_{\bif}^{\J}(f)$ exactly on the compact face $\delta_j^{\J}$. Let $\RR_+^{\J}$ be the first quadrant of $\RR^{\J}$. Then we set 
\begin{eqnarray}
I_1^{\J}&:=&\{ 1 \leq j \leq \nu (\J) \ | \ u_j^{\J} \in \Int(\RR_+^{\J}) \},
\\ 
I_2^{\J}&:=&\{ 1 \leq j \leq \nu (\J) \ | \ u_j^{\J} \in \partial \RR_+^{\J} \},\\
I_3^{\J}&:=&\{ 1 \leq j \leq \nu (\J) \ | \ u_j^{\J} \in \RR^{\J} \setminus \RR_+^{\J} \}. 
\end{eqnarray}
We thus obtain a decomposition $\{1,2,\ldots,\nu(\J)\}=I_1^{\J} \sqcup I_2^{\J} \sqcup I_3^{\J}$. For $j \in I_2^{\J}$, the primitive vector $u_j^{\J}$ must lie on one of the coordinate axes of $\RR^{\J}$. Indeed, assume that $u_j^{\J} \in \partial \RR_+^{\J}$ does not lie on any coordinate axis of $\RR^{\J}$. Then by the quasi-convenience of $f-a$ the dimension of the supporting face $\delta_j^{\J}$ of $u_j^{\J}$ in $\Gamma_{\bif}^{\J}(f)$ is less than $(\sharp \J -1)$, which contradicts our assumption $\d \delta_j^{\J}= \sharp \J -1$. For $j \in I_2^{\J}$, let $\J_j$ be the subset of $\J$ such that $\RR^{\J_j}$ is the orthogonal complement of the coordinate line $\RR u_j^{\J}$ in $\RR^{\J}$. Note that for $j \in I_2^{\J}$ we have $\sharp \J_j = \sharp \J -1$ and $\Gamma_{\infty}^{\J_j}(f)= \{0\}$ ($\Longleftrightarrow NP(f-a) \cap \RR^{\J_j}=\emptyset$). Then by using the strict non-degeneracy of $f$ along the fiber $f^{-1}(a)$, it is easy to show that
\begin{equation}
\int_{T_{\J_j}} \zeta_{f-a}(t)= (1-t^{e_j^{\J}})^{(-1)^{\sharp {\J}-1} \Vol_{\ZZ}(\delta_j^{\J})}
\end{equation}
for any $j \in I_2^{\J}$. Conversely, by the quasi-convenience of $f-a$, for any non-empty $\J \subset \{1,2,\ldots,n\}$ such that $\Gamma_{\infty}^{\J}(f)=\{0\}$ ($\Longleftrightarrow NP(f-a) \cap \RR^{\J}=\emptyset$) there exist at most one subset $\J^{\prime} \subset \{1,2,\ldots,n\}$ such that $\J \subset \J^{\prime}$, $\sharp \J^{\prime}=\sharp \J +1$ and $\Gamma_{\infty}^{\J^{\prime}}(f) \supsetneq \{0\}$. If there is no such $\J^{\prime}$ for $\J$, by \cite[Chapter I, Example (3.7)]{Oka-2} we have $\zeta_{f-a,x}(t)=1$ for any $x \in T_{\J}$ and hence
\begin{equation}
\int_{T_{\J}} \zeta_{f-a}(t)=1. 
\end{equation}
Moreover, for the empty subset $\J =\emptyset \subset \{1,2,\ldots,n\}$ we have $T_{\J}=\{ 0\} \subset \CC^n$ and $\int_{T_{\J}} \zeta_{f-a}(t)=\zeta_{f-a, 0}(t) \in \CC(t)^*$ is calculated by Varchenko's formula (\cite{Varchenko}) as
\begin{equation}
\int_{T_{\emptyset}} \zeta_{f-a}(t)=\prod_{{\J} \colon \Gamma_{\infty}^{\J}(f) \supsetneq \{0\}} \left\{ \prod_{j \in I_1^{\J}}(1-t^{e_j^{\J}})^{(-1)^{\sharp {\J}-1} \Vol_{\ZZ}(\delta_j^{\J})}\right\}.
\end{equation}
Summarizing these arguments, we obtain
\begin{eqnarray}\label{TK2} 
\tl{\zeta_f^a}(t) & = & \prod_{\J\subset \{1,2,\ldots,n\}} \int_{T_{\J}} \zeta_{f-a}(t)\\
&=& \prod_{{\J} \colon \Gamma_{\infty}^{\J}(f) \supsetneq \{0\}} \left\{ \prod_{j \in I_1^{\J} \sqcup I_2^{\J}}(1-t^{e_j^{\J}})^{(-1)^{\sharp {\J}-1} \Vol_{\ZZ}(\delta_j^{\J})}\right\} \times \prod_{{\J} \colon \Gamma_{\infty}^{\J}(f) \supsetneq \{0\}} \tl{\zeta_{f,{\J}}^a}(t). 
\end{eqnarray}
Finally, by comparing \eqref{TK1} with \eqref{TK2} we obtain the desired formula: 
\begin{equation}
\zeta_f^a(t)=\prod_{{\J} \colon \Gamma_{\infty}^{\J}(f) \supsetneq \{0\}} \left\{ \prod_{j \in I_3^{\J} }(1-t^{e_j^{\J}})^{(-1)^{\sharp {\J}-1} \Vol_{\ZZ}(\delta_j^{\J})}\right\} \times \tl{\zeta_f^a}(t). 
\end{equation}
This completes the proof.
\qed
\end{proof}

\section{Monodromy at infinity of complete intersections}\label{sec:5}

In this section, we extend our results to surjective polynomial maps $f=(f_1,f_2,\ldots,f_k) \colon \CC_x^n \longtwoheadrightarrow \CC_y^k$ ($1 \leq k \leq n$) defined by polynomials $f_1,\ldots, f_k \in \CC[x_1,\ldots, x_n]$ on $\CC^n$. It is well-known that there exists a complex hypersurface $D \subset \CC_y^k$ such that the restriction $\CC_x^n \setminus f^{-1}(D) \longtwoheadrightarrow \CC_y^k \setminus D$ of $f$ is a locally trivial fibration. We assume that the $k$-th coordinate axis
\begin{equation}
A_k :=\{y_1=y_2=\cdots =y_{k-1}=0\} \simeq \CC_{y_k}
\end{equation}
intersects $D$ at finite points: $\sharp (A_k \cap D)<+\infty$. Then there exists a finite subset $B \subset A_k \simeq \CC$ and a Zariski open subset $U \subset \CC_y^k$ such that $A_k \setminus B \subset U$ and the restriction $f^{-1}(U) \longtwoheadrightarrow U$ of $f$ is a locally trivial fibration. We denote by $B_{f,k} \subset A_k \simeq \CC$ the smallest subset of $A_k$ verifying this condition and call it the bifurcation set of $f$ on $A_k$. Set
\begin{equation}
W:=\{x \in \CC^n \ |\ f_1(x)=f_2(x)=\cdots =f_{k-1}(x)=0\}.
\end{equation}
Then by our assumption the restriction $g \colon W \longrightarrow A_k$ of $f$ induces a locally trivial fibration
\begin{equation}
W \setminus g^{-1}(B_{f,k})\longtwoheadrightarrow A_k \setminus B_{f,k}.
\end{equation}

\begin{dfn}\label{dfn:5-1}
\begin{enumerate}
\item Take a sufficiently large circle $C_R=\{y_k\in \CC \simeq A_k\ |\ |y_k|=R\}$ ($R\gg 0$) in $A_k$ such that $B_{f,k} \subset\{ y_k \in \CC \simeq A_k \ |\ |y_k|<R\}$. By restricting the locally trivial fibration $g \colon W\setminus g^{-1}(B_{f,k}) \longtwoheadrightarrow A_k \setminus B_{f,k}$ to $C_R \subset A_k \setminus B_{f,k}$, we obtain the $k$-th principal geometric monodromy at infinity
\begin{equation}
\Phi_{f,k}^{\infty} \colon g^{-1}(R) \simto g^{-1}(R)
\end{equation}
of $f=(f_1,f_2,\ldots,f_k) \colon \CC_x^n \longtwoheadrightarrow \CC_y^k$. We denote the zeta function associated with $\Phi_{f,k}^{\infty}$ by $\zeta_{f,k}^{\infty}(t) \in \CC(t)^*$ and call it the $k$-th principal monodromy zeta function at infinity of $f$.
\item For a bifurcation point $b\in B_{f,k}$ of $f$ on $A_k$, take a small circle $C_{\e}(b)=\{y_k \in \CC \simeq A_k \ |\ |y_k-b|=\e\}$ ($0<\e\ll 1$) around $b$ such that $B_{f,k} \cap \{y_k \in \CC\simeq A_k\ |\ |y_k-b|\leq\e\}=\{b\}$. We denote by $\zeta_{f,k}^b(t)\in \CC(t)^*$ the zeta function associated with the $k$-th principal geometric monodromy
\begin{equation}
\Phi_{f,k}^b \colon g^{-1}(b+\e) \simto g^{-1}(b+\e)
\end{equation}
obtained by the restriction of $g\colon W \setminus g^{-1}(B_{f,k}) \longtwoheadrightarrow A_k \setminus B_{f,k}$ to $C_{\e}(b)\subset A_k \setminus B_{f,k}$. We call $\zeta_{f,k}^b(t)$ the $k$-th principal monodromy zeta function of $f$ along the fiber $g^{-1}(b)$.
\end{enumerate}
\end{dfn}

In what follows, we always assume that $f_k \in \CC[x_1,\ldots,x_n]$ satisfies the condition $(\ast)$ (see Definition \ref{dfn:3-2}). Then the Minkowski sum 
\begin{equation}
P_{\infty}(f):=NP(f_1)+\cdots +NP(f_{k-1})+\Gamma_{\infty}(f_k)
\end{equation}
is an $n$-dimensional polytope in $\RR_v^n$. For each subset ${\J} \subset \{1,2,\ldots,n\}$ of $\{1,2,\ldots,n\}$ such that $\Gamma_{\infty}^{\J}(f_k)=\Gamma_{\infty}(f_k) \cap \RR^{\J} \supsetneq \{0\}$ ($\Longleftrightarrow$ $\d \Gamma_{\infty}^{\J}(f_k)=\sharp {\J}$), we set
\begin{equation}
I({\J}):=\{ 1 \leq j \leq k-1 \ |\ NP(f_j) \cap \RR^{\J} \neq \emptyset\} \subset \{1,2,\ldots,k-1\}
\end{equation}
and $m({\J}):=\sharp I({\J})+1$. Moreover for the $(\sharp {\J})$-dimensional cone $\RR_{\geq 0}\Gamma_{\infty}^{\J}(f_k)$ in $\RR^{\J}$ denote by $\Cone_{\J}^*$ its dual cone in $(\RR^{\J})^*$.

\begin{dfn}\label{dfn:5-2}
\begin{enumerate}
\item For a polynomial $h(x)=\sum_{v \in NP(h)}a_vx^v\in \CC[x_1,\ldots,x_n]$ ($a_v\in \CC$) on $\CC^n$ and $\J \subset \{ 1,2, \ldots , n\}$, we define a polynomial $h_{\J}(x)$ on $\CC^n$ by
\begin{equation}
h_{\J}(x):=\sum_{v \in NP(h) \cap \RR^{\J}} a_vx^v. 
\end{equation}
Moreover for each $u\in (\RR^{\J})^*$ we set 
\begin{equation}
\Gamma(h_{\J};u):=\left\{v \in NP(h) \cap \RR^{\J}\ \Bigg|\ \langle u,v \rangle=\min_{w \in NP(h) \cap \RR^{\J}} \langle u,w \rangle\right\}.
\end{equation}
\item Let $\J \subset \{ 1,2, \ldots , n\}$. For $j \in I({\J}) \sqcup \{k\}$ and $u \in (\RR^{\J})^*$, we define the $u$-part $f_j^u \in \CC[x_1,\ldots,x_n]$ of $f_j$ by
\begin{equation}
f_j^u(x):=\sum_{v \in \Gamma((f_j)_{\J};u)}a_vx^v,
\end{equation}
where $f_j(x)=\sum_{v \in NP(f_j)}a_vx^v$.
\end{enumerate}
\end{dfn}

Note that by the definition of $\Gamma_{\infty}^{\J}(f_k)$ and $\Cone_{\J}^*$ for any ${\J} \subset \{1,2,\ldots,n\}$ such that $\Gamma_{\infty}^{\J}(f_k) \supsetneq \{0\}$ and $u \in (\RR^{\J})^* \setminus \Cone_{\J}^*$ we have
\begin{equation}
\Gamma((f_k)_{\J};u)=\left\{v\in \Gamma_{\infty}^{\J}(f_k)\ \Bigg|\ \langle u,v\rangle =\min_{w \in \Gamma_{\infty}^{\J}(f_k)} \langle u, w \rangle\right\}.
\end{equation}

\begin{dfn}\label{dfn:5-3}
We say that $f=(f_1,\ldots,f_k)$ is non-degenerate at infinity if for any ${\J} \subset \{1,2,\ldots,n\}$ such that $\Gamma_{\infty}^{\J}(f_k) \supsetneq \{0\}$ and $u \in (\RR^{\J})^* \setminus \Cone_{\J}^*$ the following two subvarieties in $(\CC^*)^n$ are non-degenerate complete intersections.
\begin{gather}
\{x\in (\CC^*)^n \ |\ \text{$f_j^u(x)=0$ for any $j \in I({\J})$}\},\\
\{x\in (\CC^*)^n \ |\ \text{$f_j^u(x)=0$ for any $j \in I({\J})\sqcup \{k\}$}\}.
\end{gather}
\end{dfn}

For each subset ${\J} \subset \{1,2,\ldots,n\}$ such that $\Gamma_{\infty}^{\J}(f_k) \supsetneq \{0\}$ ($\Longleftrightarrow$ $\d \Gamma_{\infty}^{\J}(f_k)=\sharp {\J}$), consider the Minkowski sum
\begin{equation}
P_{\infty}^{\J}(f):=\sum_{j \in I({\J})} (NP(f_j)\cap \RR^{\J}) +\Gamma_{\infty}^{\J}(f_k)
\end{equation}
in $\RR^{\J}$. Then $P_{\infty}^{\J}(f)$ is a $(\sharp {\J})$-dimensional polytope in $\RR^{\J}$. Let $\gamma_1^{\J},\gamma_2^{\J},\ldots, \gamma_{n({\J})}^{\J}$ be the facets ($(\sharp {\J}-1)$-dimensional faces) of $P_{\infty}^{\J}(f)$ whose inner conormal vectors $u\neq 0\in (\RR^{\J})^*$ are contained in $(\RR^{\J})^* \setminus \Cone_{\J}^*$. For $1\leq i \leq n({\J})$, we denote by $u_i^{\J} \in (\RR^{\J})^* \setminus \Cone_{\J}^*$ the unique primitive vector in $\ZZ^{\J} \subset (\RR^{\J})^*$ which takes its minimum in $P_{\infty}^{\J}(f)$ exactly on $\gamma_i^{\J}$. For $j \in I({\J}) \sqcup\{k\}$ and $1 \leq i \leq n({\J})$, we set
\begin{equation}
\gamma(f_j)_i^{\J} :=\Gamma((f_j)_{\J};u_i^{\J}).
\end{equation}
Note that we have 
\begin{equation}
\gamma_i^{\J}=\sum_{j \in I(\J) \sqcup \{ k\} } \gamma(f_j)_i^{\J}.
\end{equation} 
For ${\J} \subset \{1,2,\ldots,n\}$ such that $\Gamma_{\infty}^{\J}(f_k) \supsetneq \{0\}$, $m(\J) \leq \sharp \J$ and $1 \leq i \leq n({\J})$, we define a positive integer $d_i^{\J}$ by
\begin{equation}
d_i^{\J}:=-\min_{w \in \Gamma_{\infty}^{\J}(f_k)} \langle u_i^{\J},w \rangle \in \ZZ_{>0}
\end{equation}
and set
\begin{equation}
K_i^{\J}:=\hspace*{-5mm}\sum_{\begin{subarray}{c}\alpha_1+\cdots +\alpha_{m({\J})} =\sharp {\J}-1\\ \text{$\alpha_q \geq 1$ for $q \leq m({\J})-1$, $\alpha_{m({\J})} \geq 0$} \end{subarray}} \hspace{-5mm}\Vol_{\ZZ}(\underbrace{\gamma(f_{j_1})_i^{\J}, \ldots, \gamma(f_{j_1})_i^{\J}}_{\text{$\alpha_1$-times}}, \ldots, \underbrace{\gamma(f_{j_{m({\J})}})_i^{\J}, \ldots, \gamma(f_{j_{m({\J})}})_i^{\J}}_{\text{$\alpha_{m({\J})}$-times}}).
\end{equation}
Here we set $I({\J}) \sqcup\{k\}=\{j_1,j_2,\ldots, j_{m(\J )-1}, k=j_{m({\J})}\}$ and
\begin{equation}
\Vol_{\ZZ}(\underbrace{\gamma(f_{j_1})_i^{\J},\ldots, \gamma(f_{j_1})_i^{\J}}_{\text{$\alpha_1$-times}}, \ldots, \underbrace{\gamma(f_{j_{m({\J})}})_i^{\J}, \ldots, \gamma(f_{j_{m({\J})}})_i^{\J}}_{\text{$\alpha_{m({\J})}$-times}})
\end{equation}
is the normalized $(\sharp {\J}-1)$-dimensional mixed volume of
\begin{equation}
\underbrace{\gamma(f_{j_1})_i^{\J},\ldots, \gamma(f_{j_1})_i^{\J}}_{\text{$\alpha_1$-times}}, \ldots, \underbrace{\gamma(f_{j_{m({\J})}})_i^{\J}, \ldots, \gamma(f_{j_{m({\J})}})_i^{\J}}_{\text{$\alpha_{m({\J})}$-times}}
\end{equation}
(see Remark \ref{rem:2-13}) with respect to the lattice $\ZZ^n \cap \LL(\gamma_i^{\J})$.

\begin{rem}\label{rem:5-4-0}
If $\sharp \J -1=0$ ($\Longrightarrow m(\J)=1$), we set
\begin{equation}
K_i^{\J}=\Vol_{\ZZ}(\underbrace{\gamma(f_k)_i^{\J}, \ldots, \gamma(f_k)_i^{\J}}_{\text{$0$ times}}):=1
\end{equation}
(in this case $\gamma(f_k)_i^{\J}$ is a point).
\end{rem}

\begin{thm}\label{thm:5-4}
Assume that $f=(f_1,\ldots,f_k)$ is non-degenerate at infinity. Then the $k$-th principal monodromy zeta function $\zeta_{f,k}^{\infty}(t)$ at infinity of $f$ is given by
\begin{equation}
\zeta_{f,k}^{\infty}(t)=\prod_{{\J} \colon \Gamma_{\infty}^{\J}(f_k) \supsetneq \{0\}, \ m({\J}) \leq \sharp {\J}} \zeta_{f,k,{\J}}^{\infty}(t),
\end{equation}
where for each ${\J} \subset \{1,2,\ldots,n\}$ such that $\Gamma_{\infty}^{\J}(f_k) \supsetneq \{0\}$ and $m({\J}) \leq \sharp {\J}$ we set
\begin{equation}
\zeta_{f,k,{\J}}^{\infty}(t):=\prod_{i=1}^{n({\J})} (1-t^{d_i^{\J}})^{(-1)^{\sharp {\J} -m({\J})}K_i^{\J}}.
\end{equation}
In particular, the Euler characteristic of the general fiber of $f \colon \CC^n \longtwoheadrightarrow \CC^k$ is equal to
\begin{equation}
\sum_{{\J} \colon \Gamma_{\infty}^{\J}(f_k) \supsetneq\{0\}, \ m({\J}) \leq \sharp {\J}} (-1)^{\sharp {\J} -m({\J})} \sum_{i=1}^{n({\J})} d_i^{\J}\cdot K_i^{\J}.
\end{equation}
\end{thm}

\begin{proof}
As in the proof of Theorem \ref{thm:3-5} (i), $\zeta_{f,k}^{\infty}(t)$ is equal to the monodromy zeta function of the nearby cycle $\psi_h(j_!(R(f_k)_!\CC_W)\in \Dbc(\{\infty\})$, where $j \colon \CC \simeq A_k \longhookrightarrow \PP^1$ is the inclusion and $h$ is a local coordinate of $\PP^1$ in a neighborhood of $\infty \in \PP^1$ such that $\infty=\{h=0\}$. Then by the standard decomposition $\CC^n =\bigsqcup_{{\J} \subset \{1,2,\ldots,n\}}T_{\J}$ of $\CC^n$ in the proof of Theorem \ref{thm:3-5} (i), we obtain a decomposition
\begin{equation}
\zeta_{f,k}^{\infty}(t)=\prod_{{\J} \subset \{1,2,\ldots,n\}} \zeta_{h, \infty}(j_!(R(f_{k, \J} )_!\CC_{W\cap T_{\J}}))(t), 
\end{equation}
where $f_{k, \J}\colon T_{\J} \longrightarrow \CC$ is the restriction of $f_k$ to $T_{\J}\simeq (\CC^*)^{\sharp \J}\subset\CC^n$. In this situation, it suffices to prove that 
\begin{equation}\label{eq:5-25}
\zeta_{h,\infty}(j_!(R(f_{k,\J})_!\CC_{W\cap T_{\J}}))(t)
=
\begin{cases}
\zeta_{f,k,{\J}}^{\infty}(t) & (\Gamma_{\infty}^{\J}(f_k) \supsetneq \{0\}, \ m({\J}) \leq \sharp {\J}),\\
1 & (\text{otherwise} ). 
\end{cases}
\end{equation}
If $\Gamma_{\infty}^{\J}(f_k)= \{0\}$, then $f_{k, \J}$ is constant on $T_{\J}$ and $\zeta_{h, \infty}(j_!(R(f_{k, \J} )_!\CC_{W\cap T_{\J}}))(t)=1$. Hence we may assume that $\Gamma_{\infty}^{\J}(f_k) \supsetneq \{0\}$ from the first. We prove the above assertion \eqref{eq:5-25} only for the case where $\J =\J_0:=\{ 1,2, \ldots , n\}$. First, let $\Sigma_1$ be the dual fan of the $n$-dimensional polytope $P_{\infty}^{\J_0}(f)=P_{\infty}(f) \subset \RR^n_v$. Next, by subdividing $\Sigma_1$ we construct a fan $\Sigma$ in $\RR^n_u=(\RR^n_v)^*$ such that the toric variety $X_{\Sigma}$ associated with it is complete and smooth. Then just as in the proof of Theorem \ref{thm:3-5} (i), we can construct a smooth variety $\tl{X_{\Sigma}}$ and the following commutative diagram of holomorphic maps: 
\begin{equation}
\xymatrix@C=20mm{ T:=T_{\J_0} \ar@{^{(}->}[r]^{\iota} \ar[d]_{f_{k,\J_0} }& \tl{X_{\Sigma}} \ar[d]^{\rho}\\
\CC \ar@{^{(}->}[r]^j & \PP^1.}
\end{equation}
Since $\rho$ is proper, we obtain
\begin{equation}
\zeta_{h, \infty}(j_!(R(f_{k, \J_0})_!\CC_{W\cap T})(t)=\zeta_{h, \infty}(R\rho_* \iota_! \CC_{W\cap T})(t). 
\end{equation}
Then by Proposition \ref{prp:2-9} for the calculation of $\zeta_{h, \infty}(R\rho_*\iota_!\CC_{W\cap T})(t)$ it suffices to calculate the monodromy zeta function of $\psi_{h \circ \rho}( \iota_! \CC_{W\cap T})$ at each point of $(h \circ \rho)^{-1}(0)=\rho^{-1}(\infty) \subset \tl{X_{\Sigma}}$. Now, let $r_1, r_2, \ldots , r_l$ be the $1$-dimensional cones (i.e. rays) in $\Sigma$ such that $r_i \setminus \{ 0\} \subset \RR^n_u \setminus \Cone_{\J_0}^*$ and $T_i$ the ($n-1$)-dimensional $T$-orbit in $X_{\Sigma}$ which corresponds to $r_i \in \Sigma$ ($i=1,2, \ldots , l$). We denote by $u_i \in \ZZ^n \setminus \{ 0\}$ the (unique) primitive vector on the ray $r_i$. If we choose an $n$-dimensional cone $\sigma_i$ in $\Sigma$ such that $r_i \prec \sigma_i$, then in the affine open subset $\CC^n(\sigma_i) \simeq \CC^n_y\subset X_{\Sigma}$ associated with $\sigma_i$ we have 
\begin{equation}
T_i=\{ y\in \CC^n(\sigma_i) \ | \ y_1 = 0, \ y_2, y_3, \ldots , y_n \neq 0\}. 
\end{equation}
Moreover for $1 \leq j \leq k$ the meromorphic extension $\tl{f_j}$ of $f_j|_T$ to $\CC^n(\sigma_i) \simeq \CC^n_y$ has the form
\begin{equation}
\tl{f_j}(y)= \frac{1}{y_1^{m_{ij}}} \times f_j^{\sigma_i}(y), 
\end{equation}
where $f_j^{\sigma_i}(y)$ is a polynomial on $\CC^n(\sigma_i)$ and we set 
\begin{equation}
m_{ij}=-\min_{w \in NP(f_j)} \langle u_i, w \rangle \ \ \in \ZZ.
\end{equation}
Note that $m_{ik} >0$ for any $1 \leq i \leq l$. By the non-degeneracy at infinity of $f=(f_1, f_2, \ldots , f_k)$, the two subvarieties
\begin{gather}
\{y\in \CC^n(\sigma_i) \ |\  f_1^{\sigma_i}(y)= \cdots =f_{k-1}^{\sigma_i}(y)=0  \},\\
\{y\in \CC^n(\sigma_i) \ |\  f_1^{\sigma_i}(y)= \cdots =f_{k-1}^{\sigma_i}(y)=f_k^{\sigma_i}(y)=0 \}
\end{gather}
in $\CC^n(\sigma_i)$ intersect $T_i$ transversally. Let us set $s_j^{\sigma_i}:=f_j^{\sigma_i}|_{T_i} \colon T_i \longrightarrow \CC$ ($i=1,2, \ldots , l$). Then by the construction of the variety $\tl{X_{\Sigma}}$ we can easily show that
\begin{equation}
\zeta_{h, \infty}(R\rho_* \iota_!\CC_{W\cap T})(t)=\prod_{i=1}^l \left( 1-t^{m_{ik}}\right)^{\chi(Z_i)}, 
\end{equation}
where we set
\begin{equation}
Z_i:= \(\bigcap_{j \in I(\J_0)} \{s_j^{\sigma_i} =0\}\) \setminus \{s_k^{\sigma_i}=0\} \subset T_i
\end{equation}
(in the case $\J =\J_0$ we have $I(\J_0)=\{ 1,2, \ldots , k-1\}$ and $m(\J_0)=k$). By Bernstein-Khovanskii-Kushnirenko's theorem (Theorem \ref{thm:2-10}), if $u_i$ is not one of the vectors $u_1^{\J_0}, u_2^{\J_0}, \ldots , u_{n(\J_0)}^{\J_0}$, the Euler characteristic $\chi (Z_i)$ of $Z_i$ is zero. Moreover by the same theorem, if $u_i=u_p^{\J_0}$ for some $1 \leq p \leq n(\J_0)$, we have 
\begin{equation}
\chi(Z_i)=(-1)^{n-k}K_p^{\J_0}=(-1)^{\sharp \J_0-m(\J_0)} K_p^{\J_0}. 
\end{equation}
Hence we obtain the desired result
\begin{equation}
\zeta_{h, \infty}(R\rho_* \iota_!\CC_{W\cap T})(t)=\prod_{i=1}^{n({\J_0})} (1-t^{d_i^{\J_0}})^{(-1)^{\sharp {\J_0} -m({\J_0})}K_i^{\J_0}}.
\end{equation}
This completes the proof.\qed
\end{proof}

To give an explicit formula for the $k$-th principal monodromy zeta function $\zeta_{f,k}^b(t)\in \CC(t)^*$ of $f$ along the fiber $g^{-1}(b)$ for $b\in B_{f,k}$, let us consider the following rational function $\tl{\zeta_{f,k}^b}(t)\in \CC(t)^*$. Let $g^{-1}(b)=\bigsqcup_{\alpha}Z_{\alpha}$ be a stratification of $g^{-1}(b)$ such that the local monodromy zeta function $\zeta_{f_k-b}(\CC_W)(t)$ is constant on each stratum $Z_{\alpha}$. Denote the value of $\zeta_{f_k-b}(\CC_W)(t)$ on $Z_{\alpha}$ by $\zeta_{\alpha}(t) \in \CC(t)^*$. Note that by the results of \cite{Kirillov}, \cite{M-T-new2} and \cite{Oka-2}, the local monodromy zeta function $\zeta_{f_k-b}(\CC_W)(t)$ is explicitly calculated in general. Then the following function $\tl{\zeta_{f,k}^b}(t)$ does not depend on the stratification $g^{-1}(b)=\bigsqcup_{\alpha}Z_{\alpha}$ of $g^{-1}(b)$.

\begin{dfn}
We set
\begin{equation}
\tl{\zeta_{f,k}^b}(t):= \dint_{g^{-1}(b)} \zeta_{f_k-b}(\CC_W)(t)=\prod_{\alpha} \{ \zeta_{\alpha}(t) \}^{\chi(Z_{\alpha})}\in \CC(t)^*
\end{equation}
and call it the finite part of $\zeta_{f,k}^b(t)$.
\end{dfn}

\begin{thm}\label{thm:5-6}
Assume that the polynomial $f_k$ is convenient and $f=(f_1,\ldots,f_k)$ is non-degenerate at infinity. Then for any bifurcation point $b \in B_{f,k}$ of $f$ on $A_k$ we have
\begin{equation}
\zeta_{f,k}^b(t)=\tl{\zeta_{f,k}^b}(t).
\end{equation}
\end{thm}

Since the proof of Theorem \ref{thm:5-6} is similar to those of Theorem \ref{thm:3-5} (ii) and \ref{thm:5-4}, we omit it. By Theorem \ref{thm:5-6}, if the polynomial $f_k$ is convenient, there is no contribution to $\zeta_{f,k}^b(t)$ from the infinity in general. In order to treat the case where $f_k$ is not convenient, we introduce the following definition.

\begin{dfn}\label{dfn:5-9}
Let $f=(f_1, f_2, \ldots, f_k)$, $g$ etc. be as above and $c \in \CC\simeq A_k$ a complex number. Then we say that $f$ is strictly non-degenerate along the fiber $f^{-1}((0,0, \ldots, c)) \simeq g^{-1}(c)$ if for any ${\J} \subset \{1,2,\ldots,n\}$ such that $\Gamma_{\infty}^{\J}(f_k) \supsetneq \{0\}$ and any non-zero $u \in (\RR^{\J})^*$ the following two subvarieties in $(\CC^*)^n$ are non-degenerate complete intersections.
\begin{gather}
\{x\in (\CC^*)^n \ |\ \text{$f_j^u(x)=0$ for any $j \in I({\J})$}\},\\
\{x\in (\CC^*)^n \ |\ \text{$f_j^u(x)=0$ for any $j \in I({\J})$ and $(f_k-c)^u(x)=0$ }\}.
\end{gather}
\end{dfn}

In what follows, we denote the constant term of $f_k$ by $a \in \CC$. Let ${\J}$ be a subset of $\{1,2,\ldots,n\}$ such that $\Gamma_{\infty}^{\J}(f_k) \supsetneq \{0\}$ ($\Longleftrightarrow$ $\d \Gamma_{\infty}^{\J}(f_k)=\sharp {\J}$). Then it is easy to see that there exist only finitely many primitive vectors $w_1^{\J}, w_2^{\J}, \ldots, w_{\nu (\J)}^{\J} \in \Int(\Cone_{\J}^*)$ in $\ZZ^{\J} \subset (\RR^{\J})^*$ such that the dimension of the Minkowski sum
\begin{equation}
\sum_{j \in I({\J})} \delta(f_j)_i^{\J}+\delta(f_k-a)_i^{\J}
\end{equation}
is $(\sharp \J -1)$, where we set
\begin{equation}
\delta(f_j)_i^{\J} :=\Gamma((f_j)_{\J};w_i^{\J})
\end{equation}
for $j \in I({\J})$ and
\begin{equation}
\delta(f_k-a)_i^{\J} :=\Gamma((f_k-a)_{\J};w_i^{\J}).
\end{equation}
For ${\J} \subset \{1,2,\ldots,n\}$ such that $\Gamma_{\infty}^{\J}(f_k) \supsetneq \{0\}$, $m(\J) \leq \sharp \J$ and $1 \leq i \leq \nu ({\J})$, we define a positive integer $e_i^{\J}$ by
\begin{equation}
e_i^{\J}:=\min_{v \in NP(f_k-a) \cap \RR^{\J}} \langle w_i^{\J}, v \rangle\in \ZZ_{> 0}
\end{equation}
and set
\begin{equation}
L_i^{\J}:=\hspace*{-5mm}\sum_{\begin{subarray}{c} \alpha_1+\cdots +\alpha_{m({\J})} =\sharp {\J}-1\\ \text{$\alpha_q \geq 1$ for $q \leq m({\J})-1$, $\alpha_{m({\J})} \geq 0$} \end{subarray}} \hspace{-10mm}\Vol_{\ZZ}(\underbrace{\delta(f_{j_1})_i^{\J}, \ldots, \delta(f_{j_1})_i^{\J}}_{\text{$\alpha_1$-times}}, \ldots, \underbrace{\delta(f_{j_{m({\J})}}-a)_i^{\J}, \ldots, \delta(f_{j_{m({\J})}}-a)_i^{\J}}_{\text{ $\alpha_{m({\J})}$-times}}),
\end{equation}
where we set $I({\J}) \sqcup\{k\}=\{j_1,j_2,\ldots, j_{m(\J )-1}, k=j_{m({\J})}\}$ as before.

As in Section \ref{sec:4}, for each ${\J} \subset \{1,2,\ldots,n\}$ consider the algebraic torus
\begin{equation}
T_{\J}:=\{x=(x_1,\ldots,x_n) \in \CC^n \ |\ x_i=0\ (i \notin {\J}), \ x_i\neq 0\ (i \in {\J})\}\simeq (\CC^*)^{\sharp {\J}}
\end{equation}
and the standard decomposition $\CC^n =\bigsqcup_{{\J} \subset \{1,2,\ldots,n\}} T_{\J}$ of $\CC^n$. 

Let $f_{k, \J} \colon T_{\J} \simeq (\CC^*)^{\sharp {\J}} \longrightarrow \CC$ be the restriction of $f_k$ to $T_{\J} \subset \CC^n$. Then for $c \in \CC$ and a subset ${\J} \subset \{1,2,\ldots,n\}$ such that $\Gamma_{\infty}^{\J}(f_k) \supsetneq \{0\}$ ($\Longleftrightarrow$ $f_{k, \J}-c$ is not constant), by taking the Euler integral of the local monodromy zeta function $\zeta_{f_{k, \J}-c}(\CC_{W \cap T_{\J}})$ over $T_{\J} \cap W \cap \{ f_k(x)-c=0 \}$ we obtain a rational function $\tl{\zeta_{f,k,{\J}}^c}(t) \in \CC(t)^*$. Note that by the results of \cite{Kirillov}, \cite{M-T-new2} and \cite{Oka-2}, the local monodromy zeta function $\zeta_{f_{k, \J}-c} (\CC_{W \cap T_{\J}})$ can be explicitly calculated in general.

\begin{thm}\label{thm:5-12}
In the situation as above, assume moreover that $f$ is strictly non-degenerate along the fiber $f^{-1}((0,0,\ldots,0,c)) \simeq g^{-1}(c)$. Then we have
\begin{enumerate}
\item If $c$ is the constant term $a$ of $f_k$, the $k$-th principal monodromy zeta function $\zeta_{f,k}^{a}(t)$ of $f$ along the fiber $g^{-1}(a)$ is given by
\begin{equation}
\zeta_{f,k}^{a}(t)=\prod_{{\J} \colon \Gamma_{\infty}^{\J}(f_k) \supsetneq \{0\}, \ m({\J}) \leq \sharp {\J}} \zeta_{f,k,{\J}}^{a}(t),
\end{equation}
where for each ${\J} \subset \{1,2,\ldots,n\}$ such that $\Gamma_{\infty}^{\J}(f_k) \supsetneq \{0\}$ and $m({\J}) \leq \sharp {\J}$ we set
\begin{equation}
\zeta_{f,k,{\J}}^{a}(t):=\prod_{i=1}^{\nu ({\J})} (1-t^{e_i^{\J}})^{(-1)^{\sharp
 {\J} -m({\J})}L_i^{\J}} \times \tl{\zeta_{f,k,{\J}}^a}(t).
\end{equation}
\item If $c \neq a$, we have
\begin{equation}
\zeta_{f,k}^{c}(t)=\prod_{{\J} \colon \Gamma_{\infty}^{\J}(f_k) \supsetneq \{0\}, \ m({\J}) \leq \sharp {\J}} \tl{\zeta_{f,k,{\J}}^c}(t).
\end{equation}
\item If $c \neq a$ and for any ${\J} \subset \{1,2,\ldots,n\}$ such that $\Gamma_{\infty}^{\J}(f_k)\supsetneq \{0\}$ we have $m({\J})=k$ (e.g.$f_1, f_2, \ldots , f_{k-1}$ are convenient), then the $k$-th principal monodromy zeta function $\zeta_{f,k}^{c}(t)$ of $f$ along the fiber $g^{-1}(c)$ is given by
\begin{equation}
\zeta_{f,k}^{c}(t)=\tl{\zeta_{f,k}^c}(t).
\end{equation}
\end{enumerate}
\end{thm}

\begin{proof}
The proofs of the assertions (i) and (ii) are similar to those of Theorem \ref{thm:4-9} and \ref{thm:5-4}. Namely for each ${\J} \subset \{1,2,\ldots,n\}$ such that $\Gamma_{\infty}^{\J}(f_k) \supsetneq \{0\}$ we construct a good toric compactification $X_{\J}$ of $T_{\J}$ as in Theorem \ref{thm:4-9} and \ref{thm:5-4}. If $m({\J})>\sharp {\J}$, then by the strict non-degeneracy of $f$ along $f^{-1}((0,0,\ldots,0,c))$ the closure of $T_{\J} \cap W$ in $X_{\J}$ does not intersect $X_{\J} \setminus T_{\J}$. In other words, $T_{\J} \cap W$ is compact and hence a finite set. This implies that we have $\tl{\zeta_{f,k,{\J}}^a}(t)=1$ for $m({\J})>\sharp {\J}$. Since the remaining parts of the proof of (i) and (ii) are completely analogous to those of Theorem \ref{thm:4-9} and \ref{thm:5-4}, we omit the detail. Finally let us prove the assertion (iii). By the strict non-degeneracy of $f$ along $f^{-1}((0,0,\ldots,0,c))$ and the assumptions in (iii) we see that $W$ intersects $T_{\J}$ transversally for any ${\J} \subset \{1,2,\ldots,n\}$ such that $\Gamma_{\infty}^{\J}(f_k) \supsetneq \{0\}$. In this case, the monodromy zeta functions $\zeta_{f_{k,\J}-c}(\CC_{W \cap T_{\J}})(t)$ and $\zeta_{f_{k}-c}(\CC_{W})(t)$ coincide each other at any point of $T_{\J} \cap W \cap \{ f_k(x)-c=0 \}$. Then the result follows from the argument in the proof of Theorem \ref{thm:4-9} (ii). \qed
\end{proof}

Finally, we shall give a generalization of Theorem \ref{thm:4-12} to the case $k\geq 2$. For each ${\J} \subset \{1,2,\ldots,n\}$ such that $\Gamma_{\infty}^{\J}(f_k) \supsetneq \{0\}$, we set
\begin{equation}
\{ i_1 < i_2 < \cdots < i_{l({\J})}\} :=\{ 1\leq i \leq \nu({\J}) \ | \ w_i^{{\J}} \notin \RR_+^{\J}\}.
\end{equation}
Then we have the following theorem.

\begin{thm}\label{thm:5-10}
In the situation as above, let $a \in \CC$ be the constant term of $f_k$ and assume that $f_k-a$ is quasi-convenient (resp. $f_1,f_2,\ldots,f_{k-1}$ are convenient). Assume moreover that $f$ is strictly non-degenerate along the fiber $f^{-1}((0,0,\ldots,0,a)) \simeq g^{-1}(a)$. Then the $k$-th principal monodromy zeta function $\zeta_{f,k}^{a}(t)$ of $f$ along the fiber $g^{-1}(a)$ is given by
\begin{equation}
\zeta_{f,k}^{a}(t)=\prod_{{\J} \colon \Gamma_{\infty}^{\J}(f_k) \supsetneq \{0\}, \ \sharp {\J} \geq k} \left\{ \prod_{q=1}^{l({\J})} (1-t^{e_{i_q}^{\J}})^{(-1)^{\sharp {\J} -k}L_{i_q}^{\J}} \right\} \times \tl{\zeta_{f,k}^a}(t).
\end{equation}
\end{thm}

\begin{proof}
As the proof of Theorem \ref{thm:4-12}, we can rewrite the result of Theorem \ref{thm:5-12} (i) by using the finite part $\tl{\zeta_{f,k}^a}(t)$ of $\zeta_{f,k}^a(t)$. For this purpose, we used also a special case of \cite[Theorem 3.12]{M-T-new2}. Since the proof is completely analogous to that of Theorem \ref{thm:4-12}, we omit the detail. \qed
\end{proof}



\begin{thebibliography}{99}

\bibitem{A'Campo}
A'Campo, N., Le nombre de Lefschetz d'une monodromie, \textit{Indag. Math.}, \textbf{35} (1973), 113-118.

\bibitem{Broughton}
Broughton, S. A., Milnor numbers and the topology of polynomial hypersurfaces, \textit{Invent. Math.}, \textbf{92} (1988), 217-241.

\bibitem{Dimca}
Dimca, A., \textit{Sheaves in topology}, Universitext, Springer-Verlag, Berlin, 2004.

\bibitem{Fulton}
Fulton, W., \textit{Introduction to toric varieties}, Princeton University Press, 1993.

\bibitem{L-N}
Garc{\'i}a L{\'o}pez, R. and N{\'e}methi, A., On the monodromy at infinity of a polynomial map, \textit{Compositio Math.}, \textbf{100} (1996), 205-231.

\bibitem{G-H}
Griffiths, P. and Harris, J., \textit{Principles of algebraic geometry}, Wiley interscience, 1994.

\bibitem{GLM2}
Gusein-Zade, S., Luengo, I. and Melle-Hern{\'a}ndez, A., Zeta functions of germs of meromorphic functions, and the Newton diagram, \textit{Funct. Anal. Appl.}, \textbf{32} (1998), 93-99.

\bibitem{GLM1}
Gusein-Zade, S., Luengo, I. and Melle-Hern{\'a}ndez, A., On the zeta-function of a polynomial at infinity, \textit{Bull. Sci. Math.}, \textbf{124} (2000), 213-224.

\bibitem{Harris}
Harris, J., \textit{Algebraic geometry}, GTM \textbf{133}, Springer-Verlag, 1992.

\bibitem{H-T-T}
Hotta, R., Takeuchi, K. and Tanisaki, T., \textit{D-modules, perverse sheaves, and representation theory}, Birkh{\"a}user Boston, 2008.

\bibitem{K-S}
Kashiwara, M. and Schapira, P., \textit{Sheaves on manifolds}, Springer-Verlag, 1990.

\bibitem{Khovanskii}
Khovanskii, A.-G., Newton polyhedra and toroidal varieties, \textit{Funct. Anal. Appl.}, \textbf{11} (1978), 289-296.

\bibitem{Kirillov}
Kirillov, A.-N., Zeta function of the monodromy for complete intersection singularities, \textit{Journal of Mathematical Sciences}, \textbf{25} (1984), 1051-1057. 

\bibitem{Kushnirenko}
Kouchnirenko, A.-G., Poly\'edres de Newton et nombres de Milnor, \textit{Invent. Math.}, \textbf{32} (1976), 1-31.

\bibitem{L-S}
Libgober, A. and Sperber, S., On the zeta function of monodromy of a polynomial map, \textit{Compositio Math.}, \textbf{95} (1995), 287-307.

\bibitem{Looijenga}
Looijenga, E., \textit{Isolated singular points of complete intersections}, London Math. Soc. Lecture Notes \textbf{77}, Cambridge University Press, 1984.

\bibitem{Massey}
Massey, D., Hypercohomology of Milnor fibers, \textit{Topology}, \textbf{35} (1996), 969-1003.

\bibitem{M-T-new1}
Matsui, Y. and Takeuchi, K., A geometric degree formula for $A$-discriminants and Euler obstructions of toric varieties, arXiv:0807.3163.

\bibitem{M-T-new2}
Matsui, Y. and Takeuchi, K., Milnor fibers over singular toric varieties and nearby cycle sheaves, arXiv:0809.3148.

\bibitem{Milnor}
Milnor, J., \textit{Singular points of complex hypersurfaces}, Princeton University Press, 1968.

\bibitem{Oda}
Oda, T., \textit{Convex bodies and algebraic geometry. An introduction to the theory of toric varieties}, Springer-Verlag, 1988.

\bibitem{Oka-2}
Oka, M., \textit{Non-degenerate complete intersection singularity}, Hermann, Paris, 1997.

\bibitem{Schurmann}
Sch{\"u}rmann, J., \textit{Topology of singular spaces and constructible sheaves}, Birkh{\"a}user, 2003.

\bibitem{S-T-1}
Siersma, D. and Tib{\u a}r, M., Singularities at infinity and their vanishing cycles, \textit{Duke Math. J.}, \textbf{80} (1995), 771-783.

\bibitem{S-T-2}
Siersma, D. and Tib{\u a}r, M., Singularities at infinity and their vanishing cycles. II. Monodromy, \textit{Publ. Res. Inst. Math. Sci.}, \textbf{36} (2000), 659-679.

\bibitem{Takeuchi}
Takeuchi, K., Perverse sheaves and Milnor fibers over singular varieties, \textit{Adv. Stud. Pure Math.}, \textbf{46} (2007), 211-222.

\bibitem{Takeuchi-2}
Takeuchi, K., Monodromy at infinity of $A$-hypergeometric functions and toric compactifications, arXiv:0812.0652.

\bibitem{Varchenko}
Varchenko, A.-N., Zeta-function of monodromy and Newton's diagram, \textit{Invent. Math.}, \textbf{37} (1976), 253-262.

\end{thebibliography}
\end{document}